\documentclass[11pt,a4paper]{article}

\usepackage{amsmath,amssymb}
\newcommand{\e}{\varepsilon}
\newcommand{\va}{\varphi}
\newcommand{\D}{\Delta}

\newcommand{\n}{\nabla}
\newcommand{\N}{\frac{N}{2}}
\newcommand{\NN}{\frac{N}{p}}
\newcommand{\q}{q_{l}}
\newcommand{\ui}{u_{l}}
\newcommand{\g}{\int_{\mathbb{R}^{N}}}
\newcommand{\p}{\partial}

\newcommand{\R}{\mathbb{R}}
\newcommand{\h}{\hookrightarrow}

\newcommand{\de}{\delta}

\newtheorem{definition}{Definition}
\newtheorem{theorem}{Theorem}
\newtheorem{proposition}{Proposition}
\newtheorem{corollaire}{Corollary}

\newtheorem{remarka}{Remark}

\newtheorem{lem}{Lemma}
\newtheorem{lemme}{Lemma}

\title{Blow-up criterion, ill-posedness and existence of  strong solution for Korteweg system with infinite energy}
\author{Boris Haspot \thanks{Basque Center of Applied Mathematics, Bizkaia Technology Park, Building 500,
E-48160, Derio (Spain), haspot@bcamath.org }}
\date{}
\begin{document}
\maketitle
\begin{abstract}
This work is devoted to the
study of the initial boundary value problem for a general
isothermal model of capillary fluids derived by J.E Dunn and
J.Serrin (1985)  (see \cite{fDS}), which can be used as a phase transition model. We will prove the existence of strong solutions in finite time with discontinuous initial density, more precisely $\ln\rho_{0}$ is in $B^{\N}_{2,\infty}(\R^{N})$. Our analysis improves the results of \cite{fDD} and \cite{fH1}, \cite{fH2} by working in space of infinite energy.  In passing our result allow to consider initial data with discontinuous interfaces, whereas in all the literature the results of existence of strong solutions consider always initial density that are continuous.  More precisely we investigate the existence of strong solution for Korteweg's system when
we authorize  jump in the pressure across some hypersurface. We obtain also  a result of ill-posedness for Korteweg system and we derive a new blow-up criterion which is the main result of this paper. More precisely we show that if we control the vacuum (i.e $\frac{1}{\rho}\in L^{\infty}_{T}(\dot{B}^{0}_{N+\e,1}(\R^{N}))$ with $\e>0$ ) then we can extend the strong solutions in finite time. It extends substantially previous results obtained for compressible equations.
\end{abstract}
\section{Introduction}
We are concerned with compressible fluids endowed with internal
capillarity. The model we consider  originates from the XIXth
century work by Van der Waals and Korteweg \cite{VW,fK} and was
actually derived in its modern form in the 1980s using the second
gradient theory, see for instance \cite{fDS,fJL,fTN}. The first investigations begin with the Young-Laplace theory which claims that the phases are separated by a hypersurface and that the jump in the pressure across the hypersurface is proportional to the curvature of the hypersurface. The main difficulty consists in describing the location and the movement of the interfaces.\\
Another major problem is to understand whether the interface behaves as a discontinuity in the state space (sharp interface) or whether the phase boundary corresponds to a more regular transition (diffuse interface, DI).
The diffuse interface models have the advantage to consider only one set of equations in a single spatial domain (the density takes into account the different phases) which considerably simplifies the mathematical and numerical study (indeed in the case of sharp interfaces, we have to treat a problem with free boundary).\\
Another approach corresponds to determine equilibrium solutions which classically consists in the minimization of the free energy functional.
Unfortunately this minimization problem has an infinity of solutions, and many of them are physically wrong (some details are given later). In order to overcome this difficulty, Van der Waals in the XIX-th century was the first to add a term of capillarity to select the physically correct solutions, modulo the introduction of a diffuse interface. This theory is widely accepted as a thermodynamically consistent model for equilibria.\\
Korteweg-type models are based on an extended version of
nonequilibrium thermodynamics, which assumes that the energy of the
fluid not only depends on standard variables but also on the
gradient of
the density. Alternatively, another way to penalize the high density variations consists in applying a zero order but non-local operator to the density gradient ( see \cite{9Ro}, \cite{5Ro}, \cite{Rohdehdr}). For more results on non local Korteweg system, we refer also to \cite{CH,Has1,Has5,Has4}.\\
Let us now consider a fluid of density $\rho\geq 0$, velocity field $u\in\R^{N}$, we are now interested in the following
compressible capillary fluid model, which can be derived from a Cahn-Hilliard like free energy (see the
pioneering work by J.- E. Dunn and J. Serrin in \cite{fDS} and also in
\cite{fA,fC,fGP}).
The conservation of mass and of momentum write:
\begin{equation}
\begin{cases}
\begin{aligned}
&\frac{\p}{\p t}\rho+{\rm div}(\rho u)=0,\\
&\frac{\p}{\p t}(\rho u)+{\rm div}(\rho
u\otimes u)-\rm div(2\mu(\rho) D (u))-\n\big(\lambda(\rho)){\rm div}u\big)+\n P(\rho)={\rm div}K,
\end{aligned}
\end{cases}
\label{3systeme}
\end{equation}
where the Korteweg tensor read as following:
\begin{equation}
{\rm div}K
=\n\big(\rho\kappa(\rho)\D\rho+\frac{1}{2}(\kappa(\rho)+\rho\kappa^{'}(\rho))|\n\rho|^{2}\big)
-{\rm div}\big(\kappa(\rho)\n\rho\otimes\n\rho\big).
\label{divK}
\end{equation}
$\kappa$ is the coefficient of capillarity and is a regular function. The term
${\rm div}K$  allows to describe the variation of density at the interfaces between two phases, generally a mixture liquid-vapor. $P$ is a general increasing pressure term.
$D (u)=(\n u+^{t}\n u)$ being the stress tensor, $\mu$ and $\lambda$ are the two Lam\'e viscosity coefficients depending on the density $\rho$) and satisfying:
$$\mu>0\;\;\mbox{and}\;\;2\mu+N\lambda\geq0.$$
Here we want to investigate the existence of strong solution for the system (\ref{3systeme}) when
we authorize  jump in the pressure across some hypersurface. To do this, we need to prove the existence of strong solutions in critical space for the scaling of the equations with initial densities which are not continuous. It will be one of main interest of this paper with a new blow-up criterion.\\
Before entering in the heart of the subject, we now want to recall the classical energy inequalities. Let $\bar{\rho}>0$ be a constant reference density (in the sequel we will assume that $\bar{\rho}=1$),   and let $\Pi $ be defined
by:
$$\Pi(s)=s\biggl(\int^{s}_{\bar{\rho}}\frac{P(z)}{z^{2}}dz-\frac{P(\bar{\rho})}{\bar{\rho}}\biggl),$$
so that $P(s)=s\Pi^{'}(s)-\Pi(s)\, ,\,\Pi^{'}(\bar{\rho})=0$. 
Multiplying the equation of momentum conservation in the system
(\ref{3systeme}) by $u$ and integrating by parts over $\R^{N}$,
we obtain the following
estimate:
\begin{equation}
\begin{aligned}
&\int_{\R^{N}}\big(\frac{1}{2}\rho
|u|^{2}+(\Pi(\rho)-\Pi(\bar{\rho}))+\frac{1}{2\kappa\rho}|\nabla\rho|^{2}\big)(t)dx
+\frac{1}{2}\int_{0}^{t}\int_{\R^{N}}\mu(\rho)
|D(u)|^{2}dxdt\\
&\leq\int_{\R^{N}}\big(\frac{|m_{0}|^{2}}{2\rho}+(\Pi(\rho_{0})-\Pi(\bar{\rho}))
+\frac{1}{2\kappa\rho_{0}}|\nabla\rho_{0}|^{2}\big)dx,
\label{3inegaliteenergie1}
\end{aligned}
\end{equation}
with:
It follows that assuming that the initial total energy is finite:\\
$${\cal E}_{0}=\int_{\R^{N}}\big(\frac{|m_{0}|^{2}}{2\rho}+(\Pi(\rho_{0})-\Pi(\bar{\rho}))
+\frac{\kappa(\rho_{0})}{2}|\nabla\rho_{0}|^{2}\big)dx<+\infty\,,$$ then we
have the a priori following bounds when $P(\rho)=a\rho^{\gamma}$:
$$(\rho-1)\in L^{\infty}(L^{\gamma}_{2}),\;\;\mbox{and}\;\;\rho |u|^{2}\in L^{1}(0,\infty,L^{1}(\R^{N})),$$
$$\sqrt{\kappa(\rho)}\n\rho\in L^{\infty}(0,\infty,L^{2}(\R^{N}))^{N},\;\;\mbox{and}\;\;\n u\in
L^{2}(0,\infty,\R^{N})^{N^{2}}.$$
We refer to section \ref{section2} for the definition of the Orlicz spaces. We now want to recall the notion of scaling for Korteweg's system. Such an
approach is now classical for incompressible Navier-Stokes equation
and yields local well-posedness (or global well-posedness for small
data) in spaces with minimal regularity.
Let us explain precisely the scaling of Korteweg's system. We can
easily check that, if $(\rho,u)$ solves (\ref{3systeme}), then
$(\rho_{\lambda},u_{\lambda})$ solves also this system:
$$\rho_{\lambda}(t,x)=\rho(\lambda^{2}t,\lambda x)\,,\,u_{\lambda}(t,x)=\lambda u(\lambda^{2}t,\lambda x)$$
provided the pressure laws $P$ have been changed into
$\lambda^{2}P$.
\begin{definition}
We say that a functional space is critical with respect to the
scaling of the equation if the associated norm is invariant under
the transformation:
$$(\rho,u)\longrightarrow(\rho_{\lambda},u_{\lambda})$$
(up to a constant independent of $\lambda$).
\end{definition}
This suggests us to choose initial data
$(\rho_{0},u_{0})$ in spaces whose norm is invariant for
all $\lambda>0$ by
$(\rho_{0},u_{0})\longrightarrow(\rho_{0}(\lambda\cdot),\lambda
u_{0}(\lambda\cdot)).$
\subsubsection*{Strong solutions}
A natural candidate is the Besov space
$B^{N/2}_{2,\infty}\times (B^{N/2-1}_{2,\infty})^{N}$,
but since $B^{N/2}_{2,\infty}$ is not included in $L^{\infty}$, we cannot
expect to get $L^{\infty}$ control on the density when
$\rho_{0}\in B^{N/2}_{2,\infty}$.  In particular it implies some difficulties in order ton control the non linear terms as the pressure. An other candidate would be $(B^{N/2}_{2,\infty}\cap L^{\infty})\times (B^{N/2-1}_{2,\infty})^{N}$, however controlling the $L^{\infty}$ norms of density is an hard task, even if the Korteweg system allows regularizing effects on the density. Indeed it appears also quadratic terms in the gradient of the density that it is necessary to control in $L^{2}(L^{\infty})$ norm to be able to obtain estimate on $\rho$ in $L^{\infty}$.\\
This is the reason why in the literature the authors consider initial density which are in Banach spaces imbedded in $L^{\infty}$. Let us briefly mention that the existence of strong solutions for $N\geq2$ is known since the works by H. Hattori and D. Li \cite{fH1,fH2}. R. Danchin and B. Desjardins in \cite{fDD} improve this result by working in critical spaces for the scaling of the equations, more precisely the initial data $(\rho_{0},\rho_{0}u_{0})$
belong to $B^{\frac{N}{2}}_{2,1}\times B^{\frac{N}{2}1}_{2,1}$ (the fact that $B^{\frac{N}{2}}_{2,1}$ is embedded in $L^{\infty}$ plays a crucial role to control the vacuum but also for some reason slinked with the multiplier space theory). In \cite{3MK}, M. Kotschote showed the existence of strong solution for the isothermal model in bounded domain by using Dore\^a-Venni Theory and $\mathcal{H}^{\infty}$ calculus.
In \cite{fH1},  we generalize the results of \cite{fDD} in the case of non isothermal Korteweg system with physical coefficients depending on the density and the temperature. More precisely we get strong solutions with initial data belonging to the critical spaces $B^{\frac{N}{2}}_{2,1}\times B^{\frac{N}{2}-1}_{2,1}\times B^{\frac{N}{2}-2}_{2,1}$ when the physical coefficients depend only on the density.
\subsubsection*{Weak solutions}
We wish to prove existence and uniqueness results for $(NHV)$ in
functions spaces very close to energy spaces. In the non isothermal
non capillary case and $P(\rho)=a\rho^{\gamma}$, with $a>0$ and
$\gamma>1$, P-L. Lions in \cite{fL2} proved the global existence of
variational solutions $(\rho,u,\theta)$ to $(NHV)$ with $\kappa=0$
for $\gamma> \frac{N}{2}$ if $N\geq 4$, $\gamma\geq \frac{3N}{N+2}$
if $N=2,3$ and initial data $(\rho_{0},m_{0})$ such that:
$$\Pi(\rho_{0})-\Pi(\bar{\rho}),\;\;\frac{|m_{0}|^{2}}{\rho_{0}}\in
L^{1}(\R^{N}),\;\;\mbox{and}\;\;\rho_{0}\theta_{0}\in
L^{1}(\R^{N}).$$ These solutions are weak solutions in the classical
sense for the equation of mass conservation and for
the equation of the momentum. Notice that the main difficulty for proving Lions' theorem consists
in exhibiting strong compactness properties of the density $\rho$ in
$L^{p}_{loc}$ spaces required to pass to the limit in the pressure
term $P(\rho)=a\rho^{\gamma}$.\\
Let us mention that Feireisl in \cite{fF}  generalized the result to
$\gamma>\frac{N}{2}$ in establishing that we can obtain renormalized
solution without imposing that $\rho\in L^{2}_{loc}$, for this he
introduces the concept of oscillation defect measure evaluating the
lost of compactness.\\
We can finally cite the result from Bresch-Desjardins
in \cite{fBD},\cite{fBD1} where they show the existence of global weak
solution for $(NHV)$ with $\kappa=0$ and with a cold pressure.  In fact they are working with specific viscosity coefficients which verify an algebraic relation. It allows them to
get good estimate on the density by using new energy inequality and to
treat by compactness all the delicate terms as the pressure. In \cite{fMV}, Mellet and Vasseur improve the results of Bresch,Desjardins by dealing with the case of a general pressure $P(\rho)=a\rho^{\gamma}$ with $\gamma>1$.\\
\\
In the case $\kappa>0$, we remark then that the density belongs to
$L^{\infty}(0,\infty,{\dot {H}}^{1}(\R^{N}))$. Hence, in contrast to
the non capillary case one can easily pass to the limit in the
pressure term. However let us emphasize at this point that the above
a priori bounds do not provide any $L^{\infty}$ control on the
density from below or from above. Indeed, even in dimension $N=2$,
$H^{1}$ functions are not necessarily locally bounded. Thus, vacuum
patches are likely to form in the fluid in spite of the presence of capillary forces, which are expected to smooth out the density. It explains why it is so difficult to obtaining the existence of global strong solution in dimension $N=2$. Indeed in order to prove the existence of global weak solution the main difficulty consists in dealing with the quadratic terms in gradient of the density appearing in the capillary tensor. Recently D. Bresch, B. Desjardins and C-K. Lin in \cite{fBDL} got
some global weak solutions for the isotherm Korteweg model with some
specific viscosity coefficients. Indeed, they assume that
$\mu(\rho)=C\rho$ with $C>0$ and $\lambda(\rho)=0$. By choosing
these specific coefficients they obtain a gain of derivatives on the
density $\rho$ where $\rho$ belongs to $L^{2}(H^{2})$.
It is easy at that time with this kind of estimate on the density $\rho$ to get strong convergence on the term of capillarity. However a new difficulty appears
concerning the loss of information on the gradient of $u$ when vacuum existing (indeed the viscosity coefficients are degenerated). It becomes
involved to pass to the limit in the term $\rho_{n} u_{n}\otimes u_{n}$. That's why
the solutions of D. Bresch, B. Desjardins and C-K. Lin require some
specific test
functions which depend on the density $\rho$. \\
In \cite{fH2}, we improve the results of R. Danchin and
B. Desjardins in \cite{fDD} and D. Bresch, B. Desjardins an C-K. Lin in \cite{fBDL} by showing the existence of global weak solution with small initial data in the energy space
for specific choices of the capillary coefficients and with general viscosity coefficient.
Comparing with the results of \cite{fBDL}, we get global weak solutions with general test function $\va\in C^{0}_{\infty}(\R^{N})$not depending on the density $\rho$. In fact we have extracted of the structure of capillarity term a new energy inequality using fractionary derivative which allows a
gain of derivative on the density $\rho$. \\
In \cite{Hprepa}, we obtain finally the existence of global weak solution for Korteweg system in any dimension with large initial data for a specific choice on the viscosity and the capillarity coefficients.
\subsection{Mathematical results}
As we explained previously, one of the main difficulty to obtain strong solutions for Korteweg's system in very general Besov spaces consists in obtaining regularizing effects on the velocity and the density without assuming that the density belongs to $L^{\infty}$. Indeed in the general case, we are brought to consider parabolic equations for the velocity and for the density with variable coefficients  what requires a control on the $L^{\infty}$ norms of $\rho$ and on $\frac{1}{\rho}$ in  $L^{\infty}$ .  To avoid these main restrictions, we will work with a specific choice on the capillarity coefficients, namely $\kappa(\rho)=\frac{\kappa}{\rho}$ with $\kappa>0$.\\
In this case, we rewrite the capillarity tensor as follows:
$$
\begin{aligned}
K(\rho)=&\rho(\n\D(\ln\rho)+\frac{1}{2}\n(|\n\ln\rho|^{2})).
\end{aligned}
$$
We refer to the appendix in section \ref{appendix} for more details on the formal computations on the tensor $K$. We now want to consider the eulerian form of the system (\ref{3systeme}) when $\kappa(\rho)=\frac{\kappa}{\rho}$ with $\kappa>0$, we obtain then by dividing by $\rho$ the momentum equation the following system:
\begin{equation}
\begin{aligned}
\begin{cases}
&\p_{t}\ln\rho+u\cdot\n\ln\rho+{\rm div}u=0,\\
&\p_{t}u+u\cdot\n u-\frac{1}{\rho}{\rm div}(\mu(\rho)D u)-\frac{1}{\rho}
\n(\lambda(\rho){\rm div}u)+\n F(\rho)=\n\D(\ln\rho)+\frac{1}{2}\n(|\n\ln\rho|^{2}),\\
&(\ln\rho,u)_{\ t=0}=(\ln\rho_{0},u_{0}).
\end{cases}
\end{aligned}
\label{NHV}
\end{equation}
with $F(\rho)$ is such that $\frac{F^{'}(\rho)}{\rho}=P^{'}(\rho)$. In the sequel we will use the following definition.
\begin{definition}
We set:
$$q=\ln\rho.$$
\end{definition}
One can now state the main results of the paper.\\
\\
In the first  theorem we prove the existence of global strong solution for (\ref{NHV}) with {\it small} initial data when we choose specific viscosity coefficients and pressure terms. We also prove the existence of strong solution in finite time with large initial data. The interest of this result will be to consider initial data with discontinuous interfaces.\\
To do this, we will  assume that: 
$$\mu(\rho)=\bar{\mu}\rho,\; \lambda(\rho)=\bar{\lambda}\rho\;\;\mbox{or}\;\;0\;\;\;\;\mbox{and}\;\;\mbox{with}\;\; \bar{\mu},\bar{\lambda},K>0$$
It leads to the following system when $\lambda=0$:
\begin{equation}
\begin{aligned}
\begin{cases}
&\p_{t}q+u\cdot\n q+{\rm div}u=0,\\
&\p_{t}u+u\cdot\n u-\D u-
\mu\n q\cdot D(u)-\lambda\n q{\rm div}u+\n F(\rho)=\n\D q+\frac{1}{2}\n(|\n q|^{2}),\\
&(q,u)_{\ t=0}=(\ln\rho_{0},u_{0}).
\end{cases}
\end{aligned}
\label{NHV1}
\end{equation}
\begin{theorem}
\label{ftheo1}Let $N\geq2$. Assume that $\mu(\rho)=\bar{\mu}\rho$, $\lambda(\rho)=\bar{\lambda}\rho$ or $0$ and $P(\rho)=K\rho$ with $\bar{\mu},\bar{\lambda},K>0$. Furthermore  we suppose that:
$$q_{0}\in B^{\NN}_{p,\infty}\;\;\mbox{and}\;\;u_{0}\in B^{\NN-1}_{p,\infty}.$$
There exists  a time $T$ such that (\ref{NHV1}) has a unique  solution  $(q,u)$ on $(0,T)$ with:
$$
\begin{aligned}
&q\in\widetilde{L}^{\infty}_{T}(B^{\N}_{2,\infty})\cap\widetilde{L}_{T}^{1}(B^{\N+2}_{2,\infty}),\;\;\mbox{and}\;\;u\in\widetilde{L}_{T}^{\infty}(B^{\N-1}_{2,\infty})
\cap\widetilde{L}_{T}^{1}(B^{\N+1}_{2,\infty}).
\end{aligned}
$$
Furthermore it exists $\e_{0}$ such that if in addition $q_{0}\in B^{\N-1}_{2,\infty}$ and:
$$\|q_{0}\|_{B^{\N-1}_{2,\infty}\cap B^{\N}_{2,\infty}}+\|u_{0}\|_{B^{\N-1}_{2,\infty}}\leq \e_{0}.$$
then the solution $(q,u)$ is global and:
\begin{equation}
\begin{aligned}
&q\in\widetilde{L}^{\infty}(B^{\N-1}_{2,\infty}\cap B^{\N}_{2,\infty})\cap\widetilde{L}^{1}(B^{\N+1}_{2,\infty}\cap B^{\N+2}_{2,\infty}),\;\;\mbox{and}\;\;u\in\widetilde{L}^{\infty}(B^{\N-1}_{2,\infty})
\cap\widetilde{L}^{1}(B^{\N+1}_{2,\infty}).
\end{aligned}
\label{impra}
\end{equation}
\end{theorem}
\begin{remarka}
We would like to mention that we could easily extend the result of strong solution in finite time to the framework of Besov spaces constructed on general $L^{p}$ spaces when the initial data are as follows:
$$q_{0}\in B^{\NN}_{p,\infty}\;\;\mbox{and}\;\;u_{0}\in B^{\NN-1}_{p,\infty},$$ 
with $1\leq p<2N$. In the case of global strong solution with small initial data in these previous spaces, we refer to \cite{CD,aHas1} where the ideas could be directly applied.
\end{remarka}
\begin{remarka}
We would like to mention that for the existence of strong solution in finite time, we could obtain the same kind of results in non homogeneous Besov space. We now want also to emphasize on the fact that when we are working with critical Besov spaces for the initial data, we generally are not able to give an estimate of the time of existence $T$. However when we work with subcritical non homogeneous Besov space for the initial data, we are able to show that the time of existence for system 8\ref{NHV1}) verifies the following inequality:
$$T^{'}\geq \frac{C}{(1+\|u_{0}\|_{\dot{B}^{\NN-1+\e}_{p,\infty}}+\|\ln\rho_{0}\|_{\dot{B}^{\NN+\e}_{p,\infty}})^{\beta}},$$
with $1\leq p<2N$ and $C,\beta$ depending on the physical coefficients and on the dimension. To do this, it is just a matter to show how $(q_{L},u_{L})$ are small in the critical Besov spaces in function of the time and of the initial data. For more details we refer to the proof of theorem \ref{ftheo1}.
\end{remarka}
\begin{remarka}
We want to point out that in the theorem \ref{ftheo1}, we solve the system (\ref{NHV1}) and not the system (\ref{NHV}). Indeed in this theorem, we do not assume any control on $q_{0}$ in $L^{\infty}$, i.e means that we have no indication on the vacuum. By this way, it is not clear that a solution from (\ref{NHV1}) is also a solution from (\ref{NHV}). This result proves in a certain way that the good variable to consider is not the density $\rho$ but rather $\ln\rho$. Furthermore we want to emphasize on the fact that this result is the first for compressible system where we are able to work with such general critical initial data (in particular $u_{0}\in B^{\N-1}_{2,\infty}$) as it is the case for incompressible Navier Stokes equations (see \cite{CMP}). In particular it is the first time up my knowledge that a compressible system has a velocity which is not Lipschitz. Indeed it is generally mandatory to control the $L^{\infty}$ norm of the density.\\
In some sense, this result means that the Korteweg system with such viscosity and capillary coefficients is closer to incompressible Navier-Stokes equations than compressible Navier Stokes system with constant viscosity coefficients.
\end{remarka}
\begin{remarka}
The main interest of this paper is to obtain some solutions which have some jump across hypersurface. It means that it could appear some change of phase with interfaces of null thickness. Indeed it is remarkable to observe that the density is not continuous in your case whereas the interfaces are assumed diffuse.
\end{remarka}
\begin{remarka}
We now want to point out the specific choice on the different physical coefficients. In particular the pressure and the viscosity coefficients have a very particular form. Indeed under this form we can check that the variable $\ln\rho$ appears everywhere but also that there is no non-linear terms in the density such that it would be necessary to obtain $L^{\infty}$ control on the density for treating these terms. In particular, it is not possible to obtain the same type of result for compressible Navier-Stokes equations as we have no regularizing effects on the density. Indeed in this case it would be not possible to deal with the following term $\n\ln\rho\cdot\n u$ in $\widetilde{L}^{1}(B^{\N-1}_{2,\infty})$.\\
In the following theorem we will prove a result of ill-posedness when we consider general pressure. Indeed in these case we need to control the $L^{\infty}$ norm of the density to deal with the pressure term of type $P(\rho)=\rho^{\gamma}$ with $\gamma$ large enough. In the case where $\gamma$ is not so large, we could extend the theorem \ref{ftheo1} by taking advantage of the regularizing effect to estimate $\n\rho^{\gamma}$ in $\widetilde{L}^{1}(B^{\N-1}_{2,\infty})$. However if we do not assume enough regularity on the divergence of the velocity, we loss the $L^{\infty}$ control on the norm of the density, it will prove in the third theorem.
\end{remarka}
We note in the sequel $\mathbb{P}=(\D)^{-1}\n{\rm div}$ and $\mathbb{Q}=(\D)^{-1}{\rm curl}{\rm div}$.
\begin{theorem}
\label{ftheo2}Let $N\geq2$.  Let $P$ a suitable smooth function. Assume that $\mu$ and $\lambda$ verify the conditions of the theorem \ref{ftheo1}. Furthermore  we suppose that:
$$q_{0}\in B^{\N}_{2,1},\;\mathbb{P}u_{0}\in B^{\N-1}_{2,2}\;\;\mbox{and}\;\;\mathbb{Q}u_{0}\in B^{\N-1}_{2,1},$$
and that $\rho_{0}$ is bounded away from the vacuum (i.e $\exists c>0$ such that $\rho_{0}\geq c$). There exists  a time $T$ such that (\ref{NHV1}) has a unique  solution  $(q,u)$ on $(0,T)$ with:
$$
\begin{aligned}
&q\in\widetilde{L}_{T}^{\infty}(B^{\N}_{2,1})\cap\widetilde{L}_{T}^{1}(B^{\N+2}_{2,1}),\;\mathbb{P}u\in\widetilde{L}^{\infty}_{T}(B^{\N-1}_{2,2})
\cap\widetilde{L}^{1}_{T}(B^{\N+1}_{2,2}),\\
&\hspace{7cm}\mbox{and}\;\;\mathbb{Q}u\in\widetilde{L}_{T}^{\infty}(B^{\N-1}_{2,1})\cap\widetilde{L}_{T}^{1}(B^{\N+1}_{2,1}).
\end{aligned}
$$
Furthermore if in addition $q_{0}\in B^{\N-1}_{2,1}$ then it exists $\e_{0}$ such that if:
$$\|q_{0}\|_{B^{\N-1}_{2,1}\cap B^{\N}_{2,1}}+\|\mathbb{P}u_{0}\|_{B^{\N-1}_{2,2}}+\|\mathbb{Q}u_{0}\|_{B^{\N-1}_{2,1}}\leq \e_{0}.$$
then the solution $(q,u)$ is global and:
$$
\begin{aligned}
&q\in\widetilde{L}^{\infty}(B^{\N-1}_{2,1}\cap B^{\N}_{2,1})\cap\widetilde{L}^{1}(B^{\N+2}_{2,1}),\;\mathbb{P}u\in\widetilde{L}^{\infty}(B^{\N-1}_{2,2})\cap\widetilde{L}^{1}(B^{\N+1}_{2,2}),\\
&\hspace{7cm}\mbox{and}\;\;\mathbb{Q}u\in\widetilde{L}^{\infty}(B^{\N-1}_{2,1})\cap\widetilde{L}^{1}(B^{\N+1}_{2,1}).
\end{aligned}
$$
\end{theorem}
\begin{remarka}
We want to point out here that the solutions of system (\ref{NHV1}) are also solutions of system (\ref{NHV}) as $q\in L^{\infty}_{T}(L^{\infty})$ which means that we control the vacuum and the $L^{\infty}$ norm from the density $\rho$. This is then sufficient by composition theorem to prove that $(\rho,u)$ is solution from (\ref{NHV}).
\end{remarka}
We now want to prove that the condition $\mathbb{Q} u_{0}\in B^{\N-1}_{2,1}$ is optimal for system (\ref{NHV}) in the sense that if we choose a third Besov index which is different from $1$ then the problem is ill-posed (in the sense of the explosion of norm for the density).
\begin{theorem}
\label{ftheo3}
Let $N\geq2$.  Let $\mu$ and $\lambda$ the viscosity coefficients are as in theorem \ref{ftheo1} and $P(\rho)=a\rho$. Then there exists a sequel of initial data such that:
\begin{itemize}
\item the sequel $(q^{n}_{0},\mathbb{P}u^{n}_{0},\mathbb{Q}u^{n}_{0})$ is uniformly bounded in
$(B^{\N}_{2,r}\cap L^{\infty})\times B^{\N-1}_{2,2}\times  B^{\N-1}_{2,r}$ with $1<r<+\infty$,
\item we have also:
$$q^{n}_{0}\in B^{\N}_{2,1},\;\mathbb{P}u^{n}_{0}\in B^{\N-1}_{2,1}\;\;\mbox{and}\;\;\mathbb{Q}u^{n}_{0}\in B^{\N-1}_{2,1},$$
\item $(q^{n}_{0},\mathbb{P}u^{n}_{0},\mathbb{Q}u^{n}_{0})$ tends weakly to $(q_{0},\mathbb{P}u_{0},\mathbb{Q}u_{0})$ in $B^{\N}_{2,r}\times B^{\N-1}_{2,2}\times  B^{\N-1}_{2,r}$ ,
\end{itemize}
such that the solutions $(q^{n},u^{n})$ associated to the initial data $(q_{0}^{n},u_{0}^{n})$ of theorem \ref{ftheo2} verify:
$$
\begin{aligned}
&\|q^{n}(t_{n},\cdot)\|_{L^{\infty}}\rightarrow+\infty\;\;\;\mbox{with}\;\;t_{n}\rightarrow0.
\end{aligned}
$$
and $(q^{n},u^{n})$ converges to a solution $(q,u)$ of (\ref{NHV1}).
\end{theorem}
\begin{remarka}
This result tell us that the problem is ill-posed for $(\ref{NHV})$ indeed in this case the solution $(q,u)$ of (\ref{NHV1}) is not a solution of (\ref{NHV}).
\end{remarka}
In the specific case where we consider that $\kappa=\mu^{2}$, we obtain from system (\ref{NHV1}) (when $\lambda=0$) the following simplified model by assuming that $v=u+\frac{\kappa}{\mu}\n\ln\rho$ (we refer to the appendix for more details on the computation or \cite{Hprepa})):
\begin{equation}
\begin{cases}
\begin{aligned}
&\p_{t}\rho-\frac{\kappa}{\mu}\D\rho=-{\rm div}(\rho v),\\
&\rho\p_{t}v +\rho u\cdot\n v-\rm div(\mu\rho\,\n v)+\n P(\rho)=0,
\end{aligned}
\end{cases}
\label{3systeme2}
\end{equation}
In the sequel, we will set $m=\rho v$ and $q{'}=\rho-1$. We have then the following theorem:
\begin{theorem}
Let $\kappa=\mu^{2}$. Let  $P$ be a suitably smooth function of the density and  $1\leq p <+\infty$.
Let $m_{0}\in \dot{B}^{\frac{N}{p}-1+\e}_{p,\infty}$ with $\e>0$ and $q^{'}_{0}\in \dot{B}^{\NN+\e}_{p,\infty}$ such that $\rho_{0}\geq c>0$. \\
There exists then a positive time $T$ such that system (\ref{3systeme2}) has a unique solution
$(q^{'},m)$ with $\rho$ bounded away from $0$ and:
$$q^{'}\in \widetilde{C}([0,T],\dot{B}^{\NN+\e}_{p,1}
)\cap \widetilde{L}^{1}_{T}(\dot{B}^{\NN+2+\e}_{p,1}),\;\;m\in \widetilde{C}([0,T];\dot{B}^{\frac{N}{p}-1+\e}_{p,1})\cap\widetilde{L}^{1}([0,T],
\dot{B}^{\NN+1+\e}_{p,1}).$$
We assume now that $P(\rho)=a\rho$ with $a>0$. Furthermore if $v_{0}\in L^{\infty}$, $\n \rho_{0}\in L^{\infty}$  and the initial data are in the energy space, it means:
$$\sqrt{\rho_{0}}u_{0}\in L^{2},\,\n\sqrt{\rho_{0}}\in L^{2}\;\;\;\mbox{and}\;\;\;(\rho_{0}-1)\in L^{1}_{2}.$$
we can then extend the solution beyond $(0,T)$ if:
\begin{equation}
(\frac{1}{\rho}-1)\in L^{\infty}_{T}(\dot{B}^{0}_{N+\e,1}),\;(\frac{1}{\sqrt{\rho}}-1)\in L^{\infty}_{T}(L^{1})\,\,\mbox{and}\;\;(\sqrt{\rho}-1)\in L^{\infty}_{T}(L^{1}).
\label{blow1}
\end{equation}
\label{theo4}
\end{theorem}
\begin{remarka}
The first condition in (\ref{blow1}) has to consider as a condition on the high frequencies and the two last one on the low frequencies.\\
In particular, it would be possible to replace the first one, by the following condition:
$$(\frac{1}{\rho}-1)\in L^{\infty}(L^{9+\e})\;\;\mbox{for}\;N=3\;\;\;\mbox{or}\;\;\;(\frac{1}{\rho}-1)\in L^{\infty}(L^{3+\e})\;\;\mbox{for}\;N=2.$$
Indeed by using energy inequality, Sobolev embedding and by seeing that $\n\rho=-\frac{1}{\rho^{\frac{3}{2}}}\n\sqrt{\rho}$, we can show in this case, that the first one is verified. In some sense these last conditions are more practical because they ask only integrability on the vacuum (and not informations in terms of Besov space). However in this case we lost the condition of scaling, we would like to explain this point in the following remark.
\end{remarka}
\begin{remarka}
This result is to compare with \cite{regular} where we show a new criterion of blow-up for compressible Navier Stokes equations when the viscosity coefficients are constant. More precisely we prove that we can extend strong solution if $P(\rho)$ is in $L^{\infty}_{T}(L^{N+1+\e})$. In particular we prove that if we have this gain of integrability on the pressure and so on the density, we are then able to control the vacuum. In our context here, it is exactly the opposite, the criterion of blow up is based on a control of the vacuum i.e $(\frac{1}{\rho}-1)\in L^{\infty}_{T}(B^{0}_{N+\e,1})$. It means that the structure of the viscosity coefficients plays really a crucial role in the understanding of the phenomena for compressible fluids. Indeed in our case, we can not apply the ideas of \cite{regular} because we lost  the so-called structure of \textit{effective velocity} introduced in \cite{regular}. However as in \cite{regular}, we are able to transfer the gain of integrability on $v$ (which is due to the fact that we have enough integrability on the pressure) to a gain of regularity on the density. It will be enough to conclude. We can also point out that as in \cite{regular} our blow up criterion (where it was the case for the pressure) asks one derivative less than the scaling of the equations on $(\frac{1}{\rho}-1)$.  It means that in some sense the subscaling (one derivative more) required on the initial data ($v_{0}\in L^{\infty}$ and $\n \rho_{0}\in L^{\infty}$) is transferred in the other sense on the blow-up criterion (indeed we need of one derivative less in terms of scaling in our blow-up condition). This fact is  very deep.\\
Furthermore by following the same arguments than in the proof of theorem \ref{theo4}, it would be possible to show a blow-up criterion which would involve $\frac{1}{\rho}$ in $L^{\infty}$ but the initial data would be assumed to belong only in Besov space with a very slight subscaling (i.e $m_{0}\in \dot{B}^{\frac{N}{p}-1+\e}_{p,\infty}$ with $\e>0$, $q^{'}_{0}\in \dot{B}^{\NN+\e}_{p,\infty}$ and no more except that the initial data will have also to be in energy space).
\end{remarka}
\begin{remarka}
We could prove easily that the result of strong solution of theorem \ref{theo4} could be extended to critical space in terms of scaling, i.e when we choose
$\e=0$. It will be proved in passing in the proof of theorem \ref{theo4}. For a result of global strong solution with small initial data for the system (\ref{3systeme2}) we refer to \cite{aHas1}, indeed the proof is a easy application of \cite{aHas1} excepted that we have to take in account the regularizing effects on the density.
\end{remarka}
\begin{remarka}
We could probably extend this previous result for more general pressure terms. However it would require more integrability on $(\frac{1}{\rho}-1)$.\\
It would be also probably possible to deal with more general viscosity and capillarity but there will be certainly an algebraic relation between these quantities.
\end{remarka}
\begin{remarka}
In fact it would be easy to improve the previous result concerning the existence of strong solution in finite time by working with initial data such that:
 $m_{0}\in \dot{B}^{\frac{N}{p}-1}_{p,2}$  and $q^{'}_{0}\in \dot{B}^{\NN}_{p,2}\cap L^{\infty}$ such that $\rho_{0}\geq c>0$. To do this we have just to take advantage of the first equation in (\ref{3systeme2}) which is a heat equation. We have just to use the kernel of the heat equation and write explicitely the solution $\rho$ to observe that the $L^{\infty}$ norm is preserved but also the vacuum.
\end{remarka}
This article is structured in the following way, first of all we recall
in the section \ref{section2} some definitions an theorems on Besov and Orlicz spaces. Next we will
concentrate in the section \ref{section4} on the proof of theorem \ref{ftheo1} and \ref{ftheo2}. In section \ref{section5}, we will prove the theorem \ref{ftheo3} of ill-posedness.  In section \ref{section6} we will deal with theorem \ref{theo4}. We postpone in appendix (see section \ref{appendix}) some technical computation on the capillarity tensor and some extensions on the previous results.
\section{Orlicz and Besov spaces}
\label{section2}
\subsection{Orlicz spaces}
We begin by describing the Orlicz space in which we will work:
$$L^{q}_{p}(\R^{N})=\{f\in L^{1}_{loc}(\R^{N})/f 1_{\{|f|\leq\delta\}}\in L^{p}(\R^{N}),\;\;
f 1_{\{|f|\geq\delta\}}\in L^{q}(\R^{N})\}$$
where $\delta$ is fixed, $\delta>0$.\\
First of all, it is not difficult to check that $L^{q}_{p}$ does not
depend on the choice of $\delta>0$ since $\frac{x^{p}}{x^{q}}$ is
bounded from above and from below on any interval
$[\delta_{1},\delta_{2}]$ with $0<\delta_{1}\leq\delta_{2}<+\infty$.
In particular we deduce that we have:
$$f^{\e}\in L^{\frac{q}{\e}}_{\frac{p}{\e}}(\R^{N})\;\;\;\mbox{if}\;\;f\in L^{q}_{p}(\R^{N})\;\;\mbox{and}
\;\;p,\,q\geq\e .$$ Obviously we get
$\mbox{meas}\{|f|\geq\delta\}<+\infty$ if $f\in L^{q}_{p}(\R^{N})$
and thus we have the embedding:
$$L^{q}_{p}(\R^{N})\subset L^{q_{1}}_{p_{1}}(\R^{N})\;\;\;\mbox{if}\;\;1\leq q_{1}\leq q<+\infty,\;\;
1\leq p\leq p_{1}<+\infty. $$ Next, we choose $\Psi$ a convex
function on $[0,+\infty)$ which is equal (or equivalent) to $x^{p}$
for $x$ small and to $x^{q}$ for $x$ large, then we can define the
space $L^{q}_{p}(\R^{N})$ as follows:
\begin{definition}We define then the Orlicz space
$L^{q}_{p}(\R^{N})$ as follows:\\\texttt{}
$$L^{q}_{p}(\R^{N})=\{f\in L^{1}_{loc}(\R^{N})/\Psi(f)\in
L^{1}(\R^{N})\}.$$
\end{definition}
We can check that $L^{q}_{p}(\R^{N})$ is a linear vector space. Now
we endow $L^{q}_{p}(\R^{N})$ with  a norm so that
$L^{q}_{p}(\R^{N})$ is a separable Banach space:
$$\|f\|_{L^{q}_{p}(\R^{N})}=\inf\{t>0/\;\;\Psi(\frac{f}{t})\leq 1\}.$$
We recall now some useful properties of Orlicz spaces.
\begin{proposition}
The following properties hold:
\begin{enumerate}
\item Dual space: If $p>1$ and $q>1$
then $(L^{q}_{p}(\R^{N}))^{'}=L^{q^{'}}_{p^{'}}(\R^{N})$ where
$q^{'}=\frac{q}{q-1},\,p^{'}=\frac{p}{p-1}$.
\item $L^{q}_{p}=L^{p}+L^{q}$ if $1\leq q\leq p<+\infty$.
\item Composition: Let $F$ be a continuous function on $\R$ such that $F(0)=0$, $F$
is differentiable at $0$ and $F(t)|t|^{-\theta}\rightarrow\alpha\ne
0$ at $t\rightarrow +\infty$. Then if $q\geq\theta$,
$$F(f)\in L^{\frac{q}{\theta}}_{p}(\R^{N})\;\;\mbox{if}\;\;f\in
L^{q}_{p}(\R^{N}).$$
\end{enumerate}
\end{proposition}
Now we can recall a property on the Orlicz space concerning the
inequality of energy.
\begin{proposition}
Let $P(\rho)=a\rho^{\gamma}$ with $a>0$ and $\gamma\geq 1$ then the function $\Pi(\rho)-\Pi(1)$ is in $ L^{1}(\R^{N})$ if and only
if $(\rho-1)\in L^{\gamma}_{2}.$
\end{proposition}
{\bf Proof:}
On the set $\{|\rho-\bar{\rho}|\leq\delta\}$, $\rho$ is bounded from
above, since $\gamma>1$ we thus deduce that $j_{\gamma}(\rho)$ is
equivalent to $|\rho-\bar{\rho}|^{2}$ on the set
$\{|\rho-\bar{\rho}|\leq\delta\}$.
Next on the set  $\{|\rho-\bar{\rho}|\geq\delta\}$, we observe that
for some $\nu\in (0,1)$ and $C\in (1,+\infty)$, we have:
$$\nu|\rho-\bar{\rho}|^{\gamma}\leq j_{\gamma}(\rho)\leq C|\rho-\bar{\rho}|^{\gamma}.$$
\hfill {$\Box$}
\subsection{Littlewood-Paley theory and Besov spaces}
Throughout the paper, $C$ stands for a constant whose exact meaning depends on the context. The notation $A\lesssim B$ means
that $A\leq CB$.
For all Banach space $X$, we denote by $C([0,T],X)$ the set of continuous functions on $[0,T]$ with values in $X$.
For $p\in[1,+\infty]$, the notation $L^{p}(0,T,X)$ or $L^{p}_{T}(X)$ stands for the set of measurable functions on $(0,T)$
with values in $X$ such that $t\rightarrow\|f(t)\|_{X}$ belongs to $L^{p}(0,T)$.
Littlewood-Paley decomposition  corresponds to a dyadic
decomposition  of the space in Fourier variables.
We can use for instance any $\varphi\in C^{\infty}(\R^{N})$ and  $\chi\in C^{\infty}(\R^{N})$ ,
supported respectively in
${\cal{C}}=\{\xi\in\R^{N}/\frac{3}{4}\leq|\xi|\leq\frac{8}{3}\}$ and $B(0,\frac{4}{3})$
such that:
$$\sum_{l\in\mathbb{Z}}\varphi(2^{-l}\xi)=1\,\,\,\,\mbox{if}\,\,\,\,\xi\ne 0,$$
and:
$$\chi(\xi)+\sum_{l\in\mathbb{N}}\varphi(2^{-l}\xi)=1\,\,\,\,\mbox{if}\,\,\,\,\forall \xi.$$
Denoting $h={\cal{F}}^{-1}\varphi$, we then define the dyadic
blocks by:
$$\D_{l}u=\varphi(2^{-l}D)u=2^{lN}\int_{\R^{N}}h(2^{l}y)u(x-y)dy\,\,\,\,\mbox{and}\,\,\,S_{l}u=\sum_{k\leq
l-1}\D_{k}u\,.$$ Formally, one can write that:
$$u=\sum_{k\in\mathbb{Z}}\D_{k}u\,.$$
This decomposition is called homogeneous Littlewood-Paley
decomposition. Let us observe that the above formal equality does
not hold in ${\cal{S}}^{'}(\R^{N})$ for two reasons:
\begin{enumerate}
\item The right hand-side does not necessarily converge in
${\cal{S}}^{'}(\R^{N})$.
\item Even if it does, the equality is not
always true in ${\cal{S}}^{'}(\R^{N})$ (consider the case of the polynomials).
\end{enumerate}
For the non homogeneous decomposition, we define the dyadic blocks as follows:
$$
\begin{aligned}
&\dot{\D}_{l}u=0\;\;\mbox{for}\;\;l\leq -2,\\
&\dot{\D}_{-1}u=\chi(D)u,\\
&\dot{\D}_{l}u=\varphi(2^{-l}D)u=2^{lN}\int_{\R^{N}}h(2^{l}y)u(x-y)dy,\;\;\mbox{for}\;\;l\geq 0.
\end{aligned}
$$
Formally, one can write that:
$$u=\sum_{k\in\mathbb{Z}}\D_{k}u\,.$$
This decomposition is called non homogeneous Littlewood-Paley
decomposition.
\subsection{Homogeneous and non homogeneous Besov spaces and first properties}
\begin{definition}
We denote by ${\cal S}_{h}^{'}$ the space of temperate distribution $u$ such that:
$$\lim S_{j}u_{j\rightarrow+\infty}=0\,\,\,\mbox{in}\;\;{\cal S}^{'}.$$
\end{definition}
\begin{definition}
For
$s\in\R,\,\,p\in[1,+\infty],\,\,q\in[1,+\infty],\,\,\mbox{and}\,\,u\in{\cal{S}}^{'}(\R^{N})$
we set:
$$\|u\|_{B^{s}_{p,q}}=(\sum_{l\in\mathbb{Z}}(2^{ls}\|\D_{l}u\|_{L^{p}})^{q})^{\frac{1}{q}}.$$
The homogeneous Besov space $B^{s}_{p,q}$ is the set of  distribution $u$ in ${\cal S}_{h}^{'}$ such that $\|u\|_{B^{s}_{p,q}}<+\infty$.
\end{definition}
\begin{definition}
For
$s\in\R,\,\,p\in[1,+\infty],\,\,q\in[1,+\infty],\,\,\mbox{and}\,\,u\in{\cal{S}}^{'}(\R^{N})$
we set:
$$\|u\|_{\dot{B}^{s}_{p,q}}=(\sum_{l\in\mathbb{Z}}(2^{ls}\|\dot{\D}_{l}u\|_{L^{p}})^{q})^{\frac{1}{q}}.$$
The non homogeneous Besov space $B^{s}_{p,q}$ is the set of temperate  distribution $u$  such that $\|u\|_{\dot{B}^{s}_{p,q}}<+\infty$.
\end{definition}
In the sequel we will give only properties on the homogeneous Besov spaces but the most of then can be generalize for non homogeneous Besov spaces.
\begin{remarka}The above definition is a natural generalization of the
nonhomogeneous Sobolev and H$\ddot{\mbox{o}}$lder spaces: one can show
that $B^{s}_{\infty,\infty}$ is the nonhomogeneous
H$\ddot{\mbox{o}}$lder space $C^{s}$ and that $B^{s}_{2,2}$ is
the nonhomogeneous space $H^{s}$.
\end{remarka}
\begin{proposition}
\label{derivation,interpolation}
The following properties holds:
\begin{enumerate}
\item there exists a constant universal $C$
such that:\\
$C^{-1}\|u\|_{B^{s}_{p,r}}\leq\|\n u\|_{B^{s-1}_{p,r}}\leq
C\|u\|_{B^{s}_{p,r}}.$
\item If
$p_{1}<p_{2}$ and $r_{1}\leq r_{2}$ then $B^{s}_{p_{1},r_{1}}\hookrightarrow
B^{s-N(1/p_{1}-1/p_{2})}_{p_{2},r_{2}}$.
\item $B^{s^{'}}_{p,r_{1}}\hookrightarrow B^{s}_{p,r}$ if $s^{'}> s$ or if $s=s^{'}$ and $r_{1}\leq r$.
\end{enumerate}
\label{interpolation}
\end{proposition}
Let now recall a few product laws in Besov spaces coming directly from the paradifferential calculus of J-M. Bony
(see \cite{BJM}) and rewrite on a generalized form in \cite{AP} by H. Abidi and M. Paicu (in this article the results are written
in the case of homogeneous sapces but it can easily generalize for the nonhomogeneous Besov spaces).
\begin{proposition}
\label{produit1}
We have the following laws of product:
\begin{itemize}
\item For all $s\in\R$, $(p,r)\in[1,+\infty]^{2}$ we have:
\begin{equation}
\|uv\|_{B^{s}_{p,r}}\leq
C(\|u\|_{L^{\infty}}\|v\|_{B^{s}_{p,r}}+\|v\|_{L^{\infty}}\|u\|_{B^{s}_{p,r}})\,.
\label{2.2}
\end{equation}
\item Let $(p,p_{1},p_{2},r,\lambda_{1},\lambda_{2})\in[1,+\infty]^{2}$ such that:$\frac{1}{p}\leq\frac{1}{p_{1}}+\frac{1}{p_{2}}$,
$p_{1}\leq\lambda_{2}$, $p_{2}\leq\lambda_{1}$, $\frac{1}{p}\leq\frac{1}{p_{1}}+\frac{1}{\lambda_{1}}$ and
$\frac{1}{p}\leq\frac{1}{p_{2}}+\frac{1}{\lambda_{2}}$. We have then the following inequalities:\\
if $s_{1}+s_{2}+N\inf(0,1-\frac{1}{p_{1}}-\frac{1}{p_{2}})>0$, $s_{1}+\frac{N}{\lambda_{2}}<\frac{N}{p_{1}}$ and
$s_{2}+\frac{N}{\lambda_{1}}<\frac{N}{p_{2}}$ then:
\begin{equation}
\|uv\|_{B^{s_{1}+s_{2}-N(\frac{1}{p_{1}}+\frac{1}{p_{2}}-\frac{1}{p})}_{p,r}}\lesssim\|u\|_{B^{s_{1}}_{p_{1},r}}
\|v\|_{B^{s_{2}}_{p_{2},\infty}},
\label{2.3}
\end{equation}
when $s_{1}+\frac{N}{\lambda_{2}}=\frac{N}{p_{1}}$ (resp $s_{2}+\frac{N}{\lambda_{1}}=\frac{N}{p_{2}}$) we replace
$\|u\|_{B^{s_{1}}_{p_{1},r}}\|v\|_{B^{s_{2}}_{p_{2},\infty}}$ (resp $\|v\|_{B^{s_{2}}_{p_{2},\infty}}$) by
$\|u\|_{B^{s_{1}}_{p_{1},1}}\|v\|_{B^{s_{2}}_{p_{2},r}}$ (resp $\|v\|_{B^{s_{2}}_{p_{2},\infty}\cap L^{\infty}}$),
if $s_{1}+\frac{N}{\lambda_{2}}=\frac{N}{p_{1}}$ and $s_{2}+\frac{N}{\lambda_{1}}=\frac{N}{p_{2}}$ we take $r=1$.
\\
If $s_{1}+s_{2}=0$, $s_{1}\in(\frac{N}{\lambda_{1}}-\frac{N}{p_{2}},\frac{N}{p_{1}}-\frac{N}{\lambda_{2}}]$ and
$\frac{1}{p_{1}}+\frac{1}{p_{2}}\leq 1$ then:
\begin{equation}
\|uv\|_{B^{-N(\frac{1}{p_{1}}+\frac{1}{p_{2}}-\frac{1}{p})}_{p,\infty}}\lesssim\|u\|_{B^{s_{1}}_{p_{1},1}}
\|v\|_{B^{s_{2}}_{p_{2},\infty}}.
\label{2.4}
\end{equation}
If $|s|<\NN$ for $p\geq2$ and $-\frac{N}{p^{'}}<s<\NN$ else, we have:
\begin{equation}
\|uv\|_{B^{s}_{p,r}}\leq C\|u\|_{B^{s}_{p,r}}\|v\|_{B^{\NN}_{p,\infty}\cap L^{\infty}}.
\label{2.5}
\end{equation}
\end{itemize}
\end{proposition}
\begin{remarka}
In the sequel $p$ will be either $p_{1}$ or $p_{2}$ and in this case $\frac{1}{\lambda}=\frac{1}{p_{1}}-\frac{1}{p_{2}}$
if $p_{1}\leq p_{2}$, resp $\frac{1}{\lambda}=\frac{1}{p_{2}}-\frac{1}{p_{1}}$
if $p_{2}\leq p_{1}$.
\end{remarka}
\begin{corollaire}
\label{produit2}
Let $r\in [1,+\infty]$, $1\leq p\leq p_{1}\leq +\infty$ and $s$ such that:
\begin{itemize}
\item $s\in(-\frac{N}{p_{1}},\frac{N}{p_{1}})$ if $\frac{1}{p}+\frac{1}{p_{1}}\leq 1$,
\item $s\in(-\frac{N}{p_{1}}+N(\frac{1}{p}+\frac{1}{p_{1}}-1),\frac{N}{p_{1}})$ if $\frac{1}{p}+\frac{1}{p_{1}}> 1$,
\end{itemize}
then we have if $u\in B^{s}_{p,r}$ and $v\in B^{\frac{N}{p_{1}}}_{p_{1},\infty}\cap L^{\infty}$:
$$\|uv\|_{B^{s}_{p,r}}\leq C\|u\|_{B^{s}_{p,r}}\|v\|_{B^{\frac{N}{p_{1}}}_{p_{1},\infty}\cap L^{\infty}}.$$
\end{corollaire}
The study of non stationary PDE's requires space of type $L^{\rho}(0,T,X)$ for appropriate Banach spaces $X$. In our case, we
expect $X$ to be a Besov space, so that it is natural to localize the equation through Littlewood-Payley decomposition. But, in doing so, we obtain
bounds in spaces which are not type $L^{\rho}(0,T,X)$ (except if $r=p$).
We are now going to
define the spaces of Chemin-Lerner in which we will work, which are
a refinement of the spaces
$L_{T}^{\rho}(B^{s}_{p,r})$.
$\hspace{15cm}$
\begin{definition}
Let $\rho\in[1,+\infty]$, $T\in[1,+\infty]$ and $s_{1}\in\R$. We set:
$$\|u\|_{\widetilde{L}^{\rho}_{T}(B^{s_{1}}_{p,r})}=
\big(\sum_{l\in\mathbb{Z}}2^{lrs_{1}}\|\D_{l}u(t)\|_{L^{\rho}(L^{p})}^{r}\big)^{\frac{1}{r}}\,.$$
We then define the space $\widetilde{L}^{\rho}_{T}(B^{s_{1}}_{p,r})$ as the set of temperate distribution $u$ over
$(0,T)\times\R^{N}$ such that 
$\|u\|_{\widetilde{L}^{\rho}_{T}(B^{s_{1}}_{p,r})}<+\infty$.
\end{definition}
We set $\widetilde{C}_{T}(\widetilde{B}^{s_{1}}_{p,r})=\widetilde{L}^{\infty}_{T}(\widetilde{B}^{s_{1}}_{p,r})\cap
{\cal C}([0,T],B^{s_{1}}_{p,r})$.
Let us emphasize that, according to Minkowski inequality, we have:
$$\|u\|_{\widetilde{L}^{\rho}_{T}(B^{s_{1}}_{p,r})}\leq\|u\|_{L^{\rho}_{T}(B^{s_{1}}_{p,r})}\;\;\mbox{if}\;\;r\geq\rho
,\;\;\;\|u\|_{\widetilde{L}^{\rho}_{T}(B^{s_{1}}_{p,r})}\geq\|u\|_{L^{\rho}_{T}(B^{s_{1}}_{p,r})}\;\;\mbox{if}\;\;r\leq\rho
.$$
\begin{remarka}
It is easy to generalize proposition \ref{produit1},
to $\widetilde{L}^{\rho}_{T}(B^{s_{1}}_{p,r})$ spaces. The indices $s_{1}$, $p$, $r$
behave just as in the stationary case whereas the time exponent $\rho$ behaves according to H\"older inequality.
\end{remarka}
In the sequel we will need of composition lemma in $\widetilde{L}^{\rho}_{T}(B^{s}_{p,r})$ spaces.
\begin{lemme}
\label{composition}
Let $s>0$, $(p,r)\in[1,+\infty]$ and $u\in \widetilde{L}^{\rho}_{T}(B^{s}_{p,r})\cap L^{\infty}_{T}(L^{\infty})$.
\begin{enumerate}
 \item Let $F\in W_{loc}^{[s]+2,\infty}(\R^{N})$ such that $F(0)=0$. Then $F(u)\in \widetilde{L}^{\rho}_{T}(B^{s}_{p,r})$. More precisely there exists a function $C$ depending only on $s$, $p$, $r$, $N$ and $F$ such that:
$$\|F(u)\|_{\widetilde{L}^{\rho}_{T}(B^{s}_{p,r})}\leq C(\|u\|_{L^{\infty}_{T}(L^{\infty})})\|u\|_{\widetilde{L}^{\rho}_{T}(B^{s}_{p,r})}.$$
\item Let $F\in W_{loc}^{[s]+3,\infty}(\R^{N})$ such that $F(0)=0$. Then $F(u)-F^{'}(0)u\in \widetilde{L}^{\rho}_{T}(B^{s}_{p,r})$. More precisely there exists a function $C$ depending only on $s$, $p$, $r$, $N$ and $F$ such that:
$$\|F(u)-F^{'}(0)u\|_{\widetilde{L}^{\rho}_{T}(B^{s}_{p,r})}\leq C(\|u\|_{L^{\infty}_{T}(L^{\infty})})\|u\|^{2}_{\widetilde{L}^{\rho}_{T}(B^{s}_{p,r})}.$$
\end{enumerate}
\end{lemme}
Now we give some result on the behavior of the Besov spaces via some pseudodifferential operator (see \cite{DFourier}).
\begin{definition}
Let $m\in\R$. A smooth function function $f:\R^{N}\rightarrow\R$ is said to be a ${\cal S}^{m}$ multiplier if for all muti-index $\alpha$, there exists a constant $C_{\alpha}$ such that:
$$\forall\xi\in\R^{N},\;\;|\p^{\alpha}f(\xi)|\leq C_{\alpha}(1+|\xi|)^{m-|\alpha|}.$$
\label{smoothf}
\end{definition}
\begin{proposition}
Let $m\in\R$ and $f$ be a ${\cal S}^{m}$ multiplier. Then for all $s\in\R$ and $1\leq p,r\leq+\infty$ the operator $f(D)$ is continuous from $B^{s}_{p,r}$ to $B^{s-m}_{p,r}$.
\label{singuliere}
\end{proposition}
Let us now give some estimates for the heat equation:
\begin{proposition}
\label{5chaleur} Let $s\in\R$, $(p,r)\in[1,+\infty]^{2}$ and
$1\leq\rho_{2}\leq\rho_{1}\leq+\infty$. Assume that $u_{0}\in B^{s}_{p,r}$ and $f\in\widetilde{L}^{\rho_{2}}_{T}
(B^{s-2+2/\rho_{2}}_{p,r})$.
Let u be a solution of:
$$
\begin{cases}
\begin{aligned}
&\p_{t}u-\mu\D u=f\\
&u_{t=0}=u_{0}\,.
\end{aligned}
\end{cases}
$$
Then there exists $C>0$ depending only on $N,\mu,\rho_{1}$ and
$\rho_{2}$ such that:
$$\|u\|_{\widetilde{L}^{\rho_{1}}_{T}(\widetilde{B}^{s+2/\rho_{1}}_{p,r})}\leq C\big(
 \|u_{0}\|_{B^{s}_{p,r}}+\mu^{\frac{1}{\rho_{2}}-1}\|f\|_{\widetilde{L}^{\rho_{2}}_{T}
 (B^{s-2+2/\rho_{2}}_{p,r})}\big)\,.$$
 If in addition $r$ is finite then $u$ belongs to $C([0,T],B^{s}_{p,r})$.
\end{proposition}
\section{Proof of theorems \ref{ftheo1} and \ref{ftheo2}}
\label{section4}
In this part we are interested in proving the theorems \ref{ftheo1} and \ref{ftheo2} of
existence of strong solutions in critical space for the scaling of the equations. We want to point out that the viscosity and the capillarity coefficients are chosen with a very specific form. This fact will be crucial in the sequel of the proof in order to obtain estimates on the density without
assuming a control on the vacuum or on the $L^{\infty}$ norm of the density. Indeed when $\mu(\rho)=\rho$ and $\kappa(\rho)=\frac{\kappa}{\rho}$ with $\kappa>0$ then the system has a specific structure and provides a new BD entropy (see \cite{Hprepa} for more details).\\
As a first step, we shall study the
linear part of the system (\ref{NHV1}) about  constant reference
density, that is:
$$
\begin{cases}
\p_{t}q+{\rm div}u=F,\\
\p_{t}u-a\D u-b\n{\rm
div}u-c\n\D q=G,
\end{cases}
\leqno{(N)}
$$
\subsection{Study of the linearized equation}
We want to prove a priori estimates in Chemin-Lerner spaces for system $(N)$
with the following hypotheses on $a,b,c,d$ which are constant:
$$
\begin{aligned}
&0<c_{1}\leq a<M_{1}<\infty,\;0<c_{2}\leq
a+b<M_{2}<\infty\;\;\mbox{and}\;\;\;0<c_{3}\leq c<M_{3}<\infty.
\end{aligned}
$$
This system has been studied by Danchin and Desjardins in \cite{fDD}, the
following proposition uses exactly the same type of arguments.
\begin{proposition}
\label{flinear3}  Let $1\leq  r\leq+\infty$ , $0\leq
s\leq1$, $(q_{0},u_{0})\in B_{2,r}^{\N+s}\times (B_{2,r}^{\N-1+s})^{N}$, and
$(F,G)\in\widetilde{L}^{1}_{T}(B^{\N+s}_{2,r})\times
(\widetilde{L}^{1}_{T}(B^{\N-1+s}_{2,r}))^{N}$.\\
Let $(q,u)\in(\widetilde{L}^{1}_{T}(B^{\N+s+2}_{2,r})\cap\widetilde{L}^{\infty}_{T}(B^{\N+s}_{2,r}))\times((\widetilde{L}^{1}_{T}(B^{\N+s+1}_{2,r}))
^{N}\cap(\widetilde{L}^{\infty}_{T}(B^{\N+s-1}_{2,r})^{N})$ be a solution of the
system $(N)$, then there exists a universal constant $C$ such that:
$$
\begin{aligned}
&\|(\n q,u)\|_{\widetilde{L}^{1}_{T}(B^{\N+1+s}_{2,r})\cap \widetilde{L}^{\infty}_{T}(B^{\N-1+s}_{2,r})}\leq C(\|(\n q_{0},u_{0})\|_{B^{\N+s}_{2,2}}+\|(\n
F,G)\|_{ \widetilde{L}^{1}_{T}(B^{\N-1+s}_{2,r})}).
\end{aligned}
$$
\end{proposition}
{\bf Proof:} \\
\\
We are going to show estimates on $q_{l}$ and
$u_{l}$. So we apply to the system the operator $\D_{l}$ , and we
have then:
\begin{eqnarray}
&&\p_{t}q_{l}+{\rm div}u_{l}=F_{l}\label{fl1}\\
&&\p_{t}u_{l}-{\rm div}(a\n u_{l})-\n(b\,{\rm div}u_{l})-\n \D
q_{l}=G_{l}\label{fl2}
\end{eqnarray}
Performing integrations by parts and usinf (\ref{fl1}) we have:
$$
\begin{aligned}
-\int_{\R^{N}}u_{l}\cdot\n\D
q_{l}dx&=\int_{\R^{N}}{\rm div}u_{l}\,\D
q_{l}dx,\\
&=-\int_{\R^{N}}\p_{t}q_{l}\,\D
q_{l}dx+\int_{\R^{N}}F_{l}\,\D
q_{l}dx,\\
&=\int_{\R^{N}}\p_{t}\n q_{l}\cdot \n
q_{l}dx+\int_{\R^{N}}\n F_{l}\cdot \n
q_{l}dx,\\
&=\frac{1}{2}\frac{d}{dt}\int_{\R^{N}}|\n
q_{l}|^{2}dx-\int_{\R^{N}}\n
q_{l}.\n F_{l}\,dx.
\end{aligned}
$$
Next, we take the inner product of (\ref{fl2}) with $u_{l}$
and by using the previous equality, we have then:
\begin{equation}
\begin{aligned}
&\frac{1}{2}\frac{d}{dt}\big(\|u_{l}\|_{L^{2}}^{2}+\int_{\R^{N}}|\n
q_{l}|^{2}dx\big)+\g(a|\n u_{l}|^{2}+b|{\rm
div}u_{l}|^{2})dx=\\
&\hspace{6cm}\g G_{l}.u_{l}\,dx+\g\n q_{l}.\n F_{l}\,dx\,.\\
\label{f3}
\end{aligned}
\end{equation}
In  order to recover some terms in $\D q_{l}$ we take the
inner product of the gradient of (\ref{fl1}) with $u_{l}$, the inner
product
scalar of (\ref{fl2}) with $\n q_{l}$ and we sum, we obtain then:
\begin{equation}
\begin{aligned}
\frac{d}{dt}\g\n q_{l}.u_{l}dx+\g (\D q_{l})^{2}dx
=&\g(G_{l}.\n q_{l}+|{\rm div}u_{l}|^{2}+u_{l}.\n F_{l}\\
&\hspace{2cm}-a\n u_{l}:\n ^{2}q_{l}-b\D q_{l}{\rm div}u
_{l})dx.
\label{f4}
\end{aligned}
\end{equation}
Let $\alpha>0$ small enough. We define:
\begin{equation}
k_{l}^{2}=\|u_{l}\|_{L^{2}}^{2}+\g(\bar{\kappa}
|\n\q|^{2}+2\alpha\n\q.\ui)dx\;.\label{f5.39}
\end{equation}
By using (\ref{f3}), (\ref{f4}) and the Young inequalities, we have
by summing and the fact that $\alpha$ is small enough:
\begin{equation}
\begin{aligned}
&\frac{1}{2}\frac{d}{dt}k_{l}^{2}+\frac{1}{2}\g(a|\n
u_{l}|^{2}+\alpha b|\D\q|^{2})dx
\lesssim\|G_{l}\|_{L^{2}}(\alpha\|\n\q\|_{L^{2}}+\|\ui\|_{L^{2}})\\
&\hspace{7cm}+\|\n F_{l}\|_{L^{2}}(\alpha\|\ui\|_{L^{2}}+\|\n\q\|_{L^{2}}).
\label{fl5}
\end{aligned}
\end{equation}
For small enough $\alpha$, we have according (\ref{f5.39}):
\begin{equation}
\frac{1}{2}k_{l}^{2}\leq\|u_{l}\|^{2}+\g\bar{\kappa}
|\n\q|^{2}dx\leq\frac{3}{2}k_{l}^{2}\;. \label{fl6}
\end{equation}
Hence according to (\ref{fl5}) and (\ref{fl6}):
$$
\begin{aligned}
\frac{1}{2}\frac{d}{dt}k_{l}^{2}+K2^{2l}k_{l}^{2} \leq&\,\,
C\,k_{l}\,(\|G_{l}\|_{L^{2}}+\|\n
F_{l}\|_{L^{2}}).\\
\end{aligned}
$$
By integrating with respect to the time, we obtain:
$$
\begin{aligned}
k_{l}(t)\leq &\,e^{-K2^{2l}t}k_{l}(0)+C\int_{0}^{t}
e^{-K2^{2l}(t-\tau)}(\|\n
F_{l}(\tau)\|_{L^{2}}+\|G_{l}(\tau)\|_{L^{2}})d\tau\;.
\end{aligned}
$$
After convolution inequalities imply that:
\begin{equation}
\begin{aligned}
\|k_{l}\|_{L^{\rho}([0,T])}\leq&\,C\big(2^{-\frac{2l}{\rho}}k_{l}(0)+2^{-2l(1+\frac{1}{\rho}-\frac{1}{\rho_{1}})}\|(\n
F_{l},G_{l})\|_{L^{\rho_{1}}_{T}(L^{2})}\big).
\label{f8}
\end{aligned}
\end{equation}
Moreover we have:
$$C^{-1}\,k_{l}\leq\|\n\q\|_{L^{2}}+\|\ui\|_{L^{2}}\leq C\,k_{l}.$$
Finally multiplying by $2^{(\N-1+s+\frac{2}{\rho})l}$, taking the $l^{r}$ norm and using (\ref{fl6}), we end up with:
$$
\begin{aligned}
\|(\n q,u)&\|_{L^{\rho}_{T}(B^{\N-1+s+\frac{2}{\rho}}_{2,r})}\leq\,\|(\n
F,G)\|_{\widetilde{L}^{\rho_{1}}_{T}(B^{\N-3+s+\frac{2}{\rho_{1}}}_{2,r})}
+\|(\n q_{0},u_{0})\|_{B^{\N-1+s}_{2,r}}.
\label{f27}
\end{aligned}
$$
It conclude the proof of the proposition.
\hfill {$\Box$}
We now want extend the result of the proposition \ref{flinear3} to the case where we include the pressure term inside of the linearized system. This point will be crucial in order to deal with the existence of global strong solution with small initial data. Indeed in this case it is very important to take in account the behavior in low frequencies and so the pressure term which is local. More precely we will consider the following linear system:
$$
\begin{cases}
\p_{t}q+{\rm div}u=F,\\
\p_{t}u-a\D u-b\n{\rm
div}u-c\n\D q+d\n q=G,
\end{cases}
\leqno{(N1)}
$$
We now want to prove a priori estimates in Chemin-Lerner spaces for system $(N1)$
with the following hypotheses on $a,b,c,d$ which are constant:
$$
\begin{aligned}
&0<c_{1}\leq a<M_{1}<\infty,\;0<c_{2}\leq
a+b<M_{2}<\infty,\;0<c_{3}\leq c<M_{3}<\infty\\
&\hspace{9cm}\;\;\mbox{and}\;\;0<c_{4}\leq d<M_{4}<\infty.
\end{aligned}
$$
This system has been studied by Danchin and Desjardins in \cite{fDD}, the
following proposition uses exactly the same type of arguments.
\begin{proposition}
\label{1flinear3}  Let $1\leq  r\leq+\infty$ , $0\leq
s\leq1$, $(q_{0},u_{0})\in (B_{2,r}^{\N-1+s}\cap B_{2,r}^{\N+s})\times (B_{2,r}^{\N-1+s})^{N}$, and
$(F,G)\in\widetilde{L}^{1}_{T}(B^{\N-1+s}_{2,r}\cap B^{\N+s}_{2,r})\times
(\widetilde{L}^{1}_{T}(B^{\N-1+s}_{2,r}))^{N}$.\\
Let $(q,u)\in\big(\widetilde{L}^{1}_{T}(B^{\N+s+1}_{2,r}\cap B^{\N+s+2}_{2,r})\cap\widetilde{L}^{\infty}_{T}(B^{\N-1+s}_{2,r}\cap B^{\N+s}_{2,r})\big)\times((\widetilde{L}^{1}_{T}(B^{\N+s+1}_{2,r}))
^{N}\cap(\widetilde{L}^{\infty}_{T}(B^{\N+s-1}_{2,r})^{N})$ be a solution of the
system $(N1)$, then there exists a universal constant $C$ such that:
$$
\begin{aligned}
&\|(\n q,q,u)\|_{\widetilde{L}^{1}_{T}(B^{\N+1+s}_{2,r})\cap \widetilde{L}^{\infty}_{T}(B^{\N-1+s}_{2,r})}\leq C(\|(\n q_{0},q_{0},u_{0})\|_{B^{\N+s}_{2,2}}+\|(\n
F,F,G)\|_{ \widetilde{L}^{1}_{T}(B^{\N-1+s}_{2,r})}).
\end{aligned}
$$
\end{proposition}
{\bf Proof:} \\
\\
It suffices to follow exactly the same lines as the proof of proposition \ref{flinear3}.
\subsection{Proof of the theorem \ref{ftheo1}:}
We now are going to prove the existence of strong solutions in critical space for system (\ref{NHV1}). In particular we recall that the main interest of theorem \ref{ftheo1} is to allow discontinuous initial data for the density, such that we can authorize discontinuous interfaces.
\subsubsection*{Existence of solutions}
We use a standard scheme:
\begin{enumerate}
\item We will use a classical iterative scheme to constructed a sequence of approximated solutions $(q^{n},u^{n})$ on a bounded interval $[0,T]$ which depend not on $n$. We will get uniform estimates on $(q^{n},u^{n})$ in:
$$E_{T}=\big(\widetilde{C}_{T}(B^{\N}_{2,\infty})\cap \widetilde{L}^{1}_{T}(
B^{\N+2}_{2,\infty})\big)\times\big(\widetilde{C}_{T}(B^{\N-1}_{2,\infty})\cap\widetilde{L}^{1}_{T}(
B^{\N+1}_{2,\infty})\big).$$
\item We will prove that the sequence $(q^{n},u^{n})$ is of Cauchy and converges to a solution of (\ref{NHV}).
\end{enumerate}
\subsubsection*{First step}
We smooth out the data as follows:
$$q_{0}^{n}=S_{n}q_{0},\;\;u_{0}^{n}=S_{n}u_{0}\;\;\;\mbox{and}\;\;\;f^{n}=S_{n}f.$$
Note that we have:
$$\forall l\in\mathbb{Z},\;\;\|\D_{l}q^{n}_{0}\|_{L^{2}}\leq\|\D_{l}q_{0}\|_{L^{2}}\;\;\;\mbox{and}\;\;\;\|q^{n}_{0}\|
_{B^{\frac{N}{2}}_{2,\infty}}\leq \|q_{0}\|_{B^{\frac{N}{2}}_{2,\infty}},$$
and similar properties for $u_{0}^{n}$ and $f^{n}$, a fact which will be used repeatedly during the next steps. Now, according \cite{fH1}, one can solve (\ref{NHV1}) with the smooth data $(q_{0}^{n},u_{0}^{n},f^{n})$.
We get a solution $(q^{n},u^{n})$ on a non trivial time interval $[0,T_{n}]$ such that:
\begin{equation}
\begin{aligned}
&q^{n}\in\widetilde{C}([0,T_{n}),B^{\N}_{2,1})\cap\widetilde{L}^{1}_{T}(
B^{\N+2}_{2,1})\;\;\mbox{and}\;\;u^{n}\in\widetilde{C}([0,T_{n}),B^{\N-1}_{2,1})\cap
\widetilde{L}^{1}_{T_{n}}
(B^{\N+1}_{2,1}).
\end{aligned}
\label{a26}
\end{equation}
\subsubsection*{Uniform bounds}
Let:
$$q^{n}=q_{L}+\bar{q}^{n},\;u^{n}=u_{L}+\bar{u}^{n},$$
where $(q_{L},u_{L})$ stands for the  solution of:
\begin{equation}
\begin{cases}
\p_{t}q_{L}+{\rm div}u_{L}=0,\\
\p_{t}u_{L}-{\cal A}u_{L}-\kappa\n(\D q_{L})=0,
\end{cases}
\label{lineaire}
\end{equation}
supplemented with initial data:
$$q_{L}(0)=q_{0}\;,\;u_{L}(0)=u_{0}.$$
Using the proposition \ref{flinear3}, we obtain the following estimates on $(q_{L},u_{L})$ for all $T>0$:
$$q_{L}\in\widetilde{C}([0,T],B^{\N}_{2,\infty})\cap\widetilde{L}^{1}_{T}(
B^{\N+2}_{2,\infty})\;\;\mbox{and}\;\;u_{L}\in\widetilde{C}([0,T],B^{\N-1}_{2,\infty})\cap
\widetilde{L}^{1}_{T}
(B^{\N+1}_{2,\infty}).$$
We let $(\bar{q}^{0},\bar{u}^{0})=(0,0)$. We now want study the behavior of $(\bar{q}_{n},\bar{u}_{n})$ where $(\bar{q}_{n},\bar{u}_{n})$ are the solution of the following system:
$$
\begin{cases}
\begin{aligned}
&\p_{t}\bar{q}^{n}+{\rm div}(\bar{u}^{n})=F_{n-1},\\
& \p_{t}\bar{u}_{n}-{\cal A}\bar{u}_{n}-\kappa\n(\D\bar{q}^{n})=G_{n-1},\\
&(\bar{q}_{n},\bar{u}_{n})_{t=0}=(0,0),
\end{aligned}
\end{cases}
\leqno{(N_{1})}
$$
where:
$$
\begin{aligned}
F_{n-1}=&-u^{n-1}\cdot\n q^{n-1},\\
=&-u_{L}\cdot\n q_{L}-\bar{u}^{n-1}\cdot\n q_{L}-u^{L}\cdot\n \bar{q}^{n-1} -\bar{u}^{n-1}\cdot\n \bar{q}^{n-1}  ,\\
G_{n-1}=&-(u^{n-1})^{*}.\n u^{n-1}+\mu\n q^{n-1}\cdot D u^{n-1}+\lambda\n q^{n-1}\,{\rm div}u^{n-1}+\frac{1}{2}\n(|\n q^{n-1}|^{2})-K\n q^{n-1}.
\end{aligned}
$$
\subsubsection*{1) First Step , Uniform Bound}
 Let $\e$ be a small
parameter and  choose $T$ small enough  so that
by using the estimate of proposition \ref{flinear3} we have:
$$
\begin{aligned}
&\|u_{L}\|_{\widetilde{L}^{1}_{T}(B^{\N+1}_{2,\infty})}+\|q_{L}\|_{\widetilde{L}^{1}_{T}(B^{\N+2}_{2,\infty})}\leq\e,\\
&\|u_{L}\|_{\widetilde{L}^{\infty}_{T}(B^{\N-1}_{2,\infty})}+\|q_{L}\|_{\widetilde{L}^{\infty}_{T}(B^{\N}_{2,\infty})}\leq
A_{0}.
\end{aligned}
\leqno{({\cal{H}}_{\e})}
$$
We are going to show by induction that:
$$\|(\bar{q}^{n},\bar{u}^{n})\|_{F_{T}}\leq\sqrt{\e}.\leqno{({\cal{P}}_{n})},$$
for $\e$ small enough with:
$$F_{T}=\big(\widetilde{C}([0,T],B^{\N}_{2,\infty})\cap\widetilde{L}^{1}_{T}(
B^{\N+2}_{2,\infty})\big)\times\big(\widetilde{C}([0,T],B^{\N-1}_{2,\infty})\cap
\widetilde{L}^{1}_{T}
(B^{\N+1}_{2,\infty})\big).$$
As $(\bar{q}^{0},\bar{u}^{0})=(0,0)$ the result is true for $n=0$. We now suppose $({\cal P}_{n-1})$ (with $n\geq 1$) true and we are going to show  $({\cal P}_{n})$.
Applying proposition \ref{flinear3}  we have:
\begin{equation}
\begin{aligned}
&\|(\bar{q}^{n},\bar{u}^{n})\|_{F_{T}}\leq C\|(\n
F_{n-1},G_{n-1})\|_{\widetilde{L}^{1}_{T}(B^{\N-1}_{2,\infty})}.
\end{aligned}
\label{fi11}
\end{equation}
Bounding the right-hand side may be
done by applying
proposition  \ref{produit1}, lemma  \ref{composition} and corollary \ref{produit2}. We begin with treating the case of $\|F_{n-1}\|_{\widetilde{L}^{1}_{T}(B^{N/2}_{2,\infty})}$, we have then:
$$\|u_{L}\cdot\n q_{L}\|_{\widetilde{L}^{1}_{T}(B^{N/2}_{2,\infty})}\leq \|u_{L}\|_{\widetilde{L}^{1}_{T}(B^{N/2+1}_{2,\infty})}\|q_{L}\|_{\widetilde{L}^{\infty}_{T}(B^{N/2}_{2,\infty})}+
\|q_{L}\|_{\widetilde{L}^{1}_{T}(B^{N/2+2}_{2,\infty})} \|u_{L}\|_{\widetilde{L}^{\infty}_{T}(B^{N/2-1}_{2,\infty})}.$$
Similarly we obtain:
$$\|u_{L}\cdot\n \bar{q}^{n-1}\|_{\widetilde{L}^{1}_{T}(B^{N/2}_{2,\infty})}\leq \|u_{L}\|_{\widetilde{L}^{1}_{T}(B^{N/2+1}_{2,\infty})}\|\bar{q}^{n-1}\|_{\widetilde{L}^{\infty}_{T}(B^{N/2}_{2,\infty})}+
\|\bar{q}^{n-1}\|_{\widetilde{L}^{\frac{4}{3}}_{T}(B^{N/2+\frac{3}{2}}_{2,\infty})} \|u_{L}\|_{\widetilde{L}^{4}_{T}(B^{N/2-\frac{1}{2}}_{2,\infty})},$$
$$\|\bar{u}^{n}\cdot\n q_{L}\|_{\widetilde{L}^{1}_{T}(B^{N/2}_{2,\infty})}\leq \|\bar{u}^{n-1}\|_{\widetilde{L}^{\frac{4}{3}}_{T}(B^{N/2+\frac{1}{2}}_{2,\infty})}\|q_{L}\|_{\widetilde{L}^{4}_{T}(B^{N/2+\frac{1}{2}}_{2,\infty})}+\|q_{L}\|_{\widetilde{L}^{1}_{T}(B^{N/2+2}_{2,\infty})} \|\bar{u}^{n-1}\|_{\widetilde{L}^{\infty}_{T}(B^{N/2-1}_{2,\infty})},$$
and:
$$\|\bar{u}^{n-1}\cdot\n \bar{q}^{n-1}\|_{\widetilde{L}^{1}_{T}(B^{N/2}_{2,\infty})}\leq \|\bar{u}^{n-1}\|_{\widetilde{L}^{1}_{T}(B^{N/2+1}_{2,\infty})}\|\bar{q}^{n-1}\|_{\widetilde{L}^{\infty}_{T}(B^{N/2}_{2,\infty})}+\|\bar{q}^{n-1}\|_{\widetilde{L}^{1}_{T}(B^{N/2+2}_{2,\infty})} \|\bar{u}^{n-1}\|_{\widetilde{L}^{\infty}_{T}(B^{N/2-1}_{2,\infty})}.$$
By using the previous inequalities and $({\cal H}_{\e})$, we obtain that:
\begin{equation}
\|F_{n}\|_{L^{1}_{T}(B^{N/2})}\leq C(2A_{0}\e+2\e^{\frac{3}{2}}+2\sqrt{\e}\e^{\frac{1}{4}}+2\e).
\label{Fn}
\end{equation}
Next we want to control  $\|G_{n}\|_{\widetilde{L}^{1}(B^{\N-1}_{2,\infty})}$. According to
propositions \ref{produit1}, corollary \ref{produit2} and \ref{flinear3}, we have:
$$
\begin{aligned}
&\|(u^{n-1})^{*}.\n
u^{n-1}\|_{L^{1}_{T}(B^{\N-1}_{2,\infty})}\lesssim\|u^{n-1}\|_{L^{\frac{4}{3}}_{T}(B^{\N+\frac{1}{2}}_{2,\infty})}\|u^{n-1}\|_{L^{4}_{T}(B^{\N-\frac{1}{2}}_{2,\infty})},\\
&\|\n(|\n q^{n-1}|^{2})\|_{\widetilde{L}^{1}_{T}(B^{\N-1}_{2,\infty})}\lesssim\||\n q^{n-1}|^{2}
\|_{\widetilde{L}^{1}_{T}(B^{\N}_{2,\infty})},\\
&\hspace{3,6cm}\lesssim \|\n q^{n-1}
\|_{\widetilde{L}^{\frac{4}{3}}_{T}(B^{\N+\frac{1}{2}}_{2,\infty})}\|\n q^{n-1}
\|_{\widetilde{L}^{4}_{T}(B^{\N-\frac{1}{2}}_{2,\infty})},\\
&\hspace{3,6cm}\lesssim \|q^{n-1}
\|_{\widetilde{L}^{\frac{4}{3}}_{T}(B^{\N+\frac{3}{2}}_{2,\infty})}\|q^{n-1}
\|_{\widetilde{L}^{4}_{T}(B^{\N+\frac{1}{2}}_{2,\infty})}.
\end{aligned}
$$
We proceed similarly for the other terms and we obtain by using (\ref{fi11}) and the different previous inequalities:
$$\|(\bar{q}_{n+1},\bar{u}_{n+1})\|_{F_{T}}\leq C\sqrt{\e}(\sqrt{\e}A_{0}+\sqrt{\e}+\e^{\frac{1}{4}}).$$
By taking $T$ and $\e$ small enough  we have $({\cal{P}}_{n+1})$, so we
have shown by induction that $(q^{n},u^{n})$ is bounded
in $F_{T}$.\\
\\
In the case where we want to obtain global strong solution with small initial data, we need to take in account the low frequencies. That is why we will include the pressure term in the linear part and we will use proposition \ref{1flinear3} to obtain estimate on $(q_{L},u_{L})$ where $(q_{L},u_{L})$ here stands for the  solution of:
\begin{equation}
\begin{cases}
\p_{t}q_{L}+{\rm div}u_{L}=0,\\
\p_{t}u_{L}-{\cal A}u_{L}+K\n q_{L}-\kappa\n(\D q_{L})=0,
\end{cases}
\label{lineaire}
\end{equation}
In particular it explains why we need of additional regularity on the initial density, i.e $q_{0}\in B^{\N-1}_{2,\infty}$. The rest of the proof for this case is similar to the case with large initial data.
\subsubsection*{Second Step: Convergence of the
sequence}
 We will show
that $(q^{n},u^{n})$ is a Cauchy sequence in the Banach
space $F_{T}$, hence converges to some
$(q,u)\in F_{T}$.\\
Let:
$$\delta q^{n}=q^{n+1}-q^{n},\;\delta u^{n}=u^{n+1}-u^{n}.$$
The system verified by $(\de q^{n},\de u^{n})$ reads:
$$
\begin{cases}
\begin{aligned}
&\p_{t}\delta q^{n}+{\rm div}\delta u^{n}=F_{n}-F_{n-1},\\
&\p_{t}\delta u^{n}-\mu\D\delta u^{n}-(\lambda+\mu)\n{\rm div}\delta u^{n} -\n \D \delta q^{n}=G_{n}-G_{n-1},\\
&\delta q^{n}(0)=0\;,\;\delta u^{n}(0)=0,
\end{aligned}
\end{cases}
$$
Applying propositions \ref{flinear3}, and using $({\cal{P}}_{n})$, we get:
$$
\begin{aligned}
\|(\de q^{n},\de u^{n},\de{\cal T}^{n})\|_{F_{T}}\leq\;&
C(\|F_{n}-F_{n-1}\|_{L^{1}_{T}(B^{N/2}_{2,\infty})}+\|G_{n}-G_{n-1}\|_{L^{1}_{T}(B^{N/2-1}_{2,\infty})}).
\end{aligned}
$$
And by the same type of estimates as before, we get:
$$\|(\de q^{n},\de u^{n})\|_{F_{T}}\leq
C\sqrt{\e}(1+A_{0})^{3}\|(\de q^{n-1},\de u^{n-1})\|_{F_{T}}.$$ So in taking $\e$ enough small we have that
$(q^{n},u^{n})$ is Cauchy sequence, so the limit
$(q,u)$ is in $F_{T}$ and we verify easily that this is a
solution of the system.
\subsubsection*{Third step: Uniqueness}
Now, we are going to prove the uniqueness of the solution in the following space:
$$\widetilde{F}^{\N}_{T}=\big(\widetilde{C}([0,T],B^{\N}_{2,\infty})\cap\widetilde{L}^{2}_{T}(
B^{\N+1}_{2,\infty})\big)\times\big(\widetilde{C}([0,T],B^{\N-1}_{2,\infty})\cap
\widetilde{L}^{2}_{T}
(B^{\N}_{2,\infty})\big).$$
Suppose that $(q_{1},u_{1})$ and $(q_{2},u_{2})$ are solutions with the same initial conditions, and
$(q_{1},u_{1})$ corresponds to the previous
solution.\\
We set then:
$$\de q=q_{2}-q_{1}\;\;\;\mbox{and}\;\;\;\de u=u_{2}-u_{1}.$$
$(\de q,\de u)$ satisfy the following system:
$$
\begin{cases}
\begin{aligned}
&\p_{t}\delta q+{\rm div}\delta
u=F_{2}-F_{1},\\
&\p_{t}\delta u-\mu\D\delta u-(\lambda+\mu)\n{\rm div}\delta u -\n \D \delta q=G_{1}-G_{2},\\
&\delta q(0)=0\;,\;\delta u(0)=0.
\end{aligned}
\end{cases}
$$
We now apply proposition \ref{flinear3} to the previous system, and by using the same type of estimates than in the part on the contraction, we
obtain:
$$
\begin{aligned}
&\|(\de q,\de u)\|_{\widetilde{F}^{\N}_{T_{1}}}\lesssim(\|q_{1}\|_{\widetilde{L}^{2}{T_{1}}(B^{\N+1}_{2,\infty})}+\|q_{2}\|_{\widetilde{L}{T_{1}}^{2}(B^{\N+1}_{2,\infty})}+
\|u_{1}\|_{\widetilde{L}^{2}{T_{1}}(B^{\N+1}_{2,\infty})}+\|u_{2}\|_{\widetilde{L}{T_{1}}^{2}(B^{\N}_{2,\infty})})\\
&\hspace{11cm}\times\|(\de
q,\de u,\de{\cal T})\|_{\widetilde{F}^{\N}_{T_{1}}}.
\end{aligned}
$$
We have then for $T_{1}$ small enough: $(\de q,\de u)=(0,0)$ on $[0,T_{1}]$ and by connectivity we finally
conclude that:
$$q_{1}=q_{2},\;u_{1}=u_{2}\;\;\mbox{on}\;\;[0,T].$$ \hfill {$\Box$}
\subsection{Proof of the theorem \ref{ftheo2}}
In this section we have to deal with general pressure. It means that we need of additional information to control the pressure in Besov space. Indeed in order to
use some composition law of paraproduct on the pressure it is necessary to control the $L^{\infty}$ norm of the density $\rho$. To do this we need of additional regularity assumption on the initial data to be able to estimate the $L^{\infty}$ norm on the density $\rho$.\\
The proof follows the same lines as the proof of theorem \ref{ftheo1} except that we obtain new estimates on $(q_{L},u_{L})$, indeed by using proposition \ref{flinear3},  $(q_{L},u_{L})$ is in $E_{T}$ with:
$$E_{T}=\big(\widetilde{C}([0,T],B^{\N}_{2,2})\cap\widetilde{L}^{1}_{T}(
B^{\N+2}_{2,2})\big)\times\big(\widetilde{C}([0,T],B^{\N-1}_{2,2})\cap
\widetilde{L}^{1}_{T}
(B^{\N+1}_{2,2})\big).$$
However in order to estimate via the mass equation the density in $L^{\infty}$, we need to control ${\rm div}u$ in $L^{1}(L^{\infty})$.That is why as ${\cal Q}u_{0}$ belongs to $B^{\N-1}_{2,1}$, we are able to prove that ${\rm div}u_{L}\in\widetilde{L}^{1}(B^{\N}_{2,1})$. Indeed applying the operator ${\rm div}$ to the momentum equation of system (\ref{lineaire}) and letting $v_{L}={\rm div}u_{L}$, we obtain the following system:
\begin{equation}
\begin{cases}
&\p_{t}c+\D v_{L}=0,\\
&\p_{t}v_{L}-(2\mu+\lambda)\D v_{L}-\kappa\D c=0,\\
&(c,v_{L})_{/ t=0}=(\D \ln\rho_{0}, {\rm div}u_{0}),
\end{cases}
\label{div1}
\end{equation}
with $c=\D \ln\rho$.  We now recall a lemma obtained in \cite{fDD}.
\begin{lem}
\label{petitlemme}
Let $s\in\R$ and $r\in[1,+\infty]$. Suppose that $v_{0}\in B^{s}_{2,r}$. Then the system (\ref{div})
has a unique solution $v_{L}$ in $\widetilde{C}(B^{s}_{2,r})\cap \widetilde{L}^{1}(B^{s+2}_{2,r})$ and we have:
$$\|v_{L}\|_{\widetilde{L}^{\infty}(B^{s}_{2,r})}+\|v_{L}\|_{\widetilde{L}^{1}(B^{s+2}_{2,r})}\leq C\|v_{0}\|_{B^{s}_{2,r}},$$
where $C$ depends on $\mu$, $\lambda$ and $\kappa$.
\end{lem}
As a consequence of lemma \ref{petitlemme}, we obtain that ${\rm div}u_{L}$ belongs to $\widetilde{L}^{1}(B^{\N}_{2,1})$. Unfortunately it is not sufficient to hope controlling the density in norm $L^{\infty}$. Indeed via the transport equation, we need also to control the divergence of the velocity in $L^{1}(L^{\infty})$. To do this, it suffices to estimate $\bar{u}^{n}$ in $\widetilde{L}^{1}(B^{\N}_{2,1})$. The idea is then to prove that $\bar{u}^{n}$ is more regular on the third index than $u_{L}$. This type of result is well-known in the case of Navier-Stokes (see \cite{CP}), we now want to adapt the spirit of this result to our case.\\
We now have to prove that $(\bar{q}^{n},\bar{u}^{n})$ is bounded in $E^{1}_{T}$ with:
$$E^{1}_{T}=\big(\widetilde{C}([0,T],B^{\N}_{2,1})\cap\widetilde{L}^{1}_{T}(
B^{\N+2}_{2,1})\big)\times\big(\widetilde{C}([0,T],B^{\N-1}_{2,1})\cap
\widetilde{L}^{1}_{T}
(B^{\N+1}_{2,1})\big).$$
It means in particular that we have a type of regularizing effects on $(\bar{q}^{n},u^{n})$ on the third index of the Besov spaces. More precisely we are going to show by induction as in theorem \ref{ftheo1} that:
$$\|(\bar{q}^{n},\bar{u}^{n})\|_{E^{1}_{T}}\leq\e.\leqno{({\cal{P}}_{n})},$$
for $\e$ small enough.\\
As $(\bar{q}^{0},\bar{u}^{0})=(0,0)$ the result is true for $n=0$.We now suppose $({\cal P}_{n-1})$ (with $n\geq 1$) true and we are going to show  $({\cal P}_{n})$.
Applying proposition \ref{flinear3}  we have:
\begin{equation}
\begin{aligned}
&\|(\bar{q}^{n},\bar{u}^{n})\|_{F_{T}}\leq C\|(\n
F_{n},G_{n})\|_{\widetilde{L}^{1}_{T}(B^{\N-1}_{2,1})}.
\end{aligned}
\label{fi1}
\end{equation}
Bounding the right-hand side may be
done by applying
proposition  \ref{produit1}, lemma  \ref{composition} and corollary \ref{produit2}.  We want only to give an example of how to control  $\|G_{n}\|_{\widetilde{L}^{1}(B^{\N-1}_{2,\infty})}$.  We will treat the case of $u^{n-1})^{*}.\n
u^{n-1}$ and $\n(|\n q^{n-1}|^{2})$. According to
propositions \ref{produit1}, corollary \ref{produit2} and \ref{flinear3}, we have:
$$
\begin{aligned}
&\|(u^{n-1})^{*}.\n
u^{n-1}\|_{L^{1}_{T}(B^{\N-1}_{2,1})}\lesssim\|u^{n-1}\|_{L^{\infty}_{T}(B^{\N-1}_{2,2})}\|u^{n-1}\|_{L^{1}_{T}(B^{\N+1}_{2,2})},\\
&\|\n(|\n q^{n-1}|^{2})\|_{\widetilde{L}^{1}_{T}(B^{\N-1}_{2,\infty})}\lesssim\||\n q^{n-1}|^{2}
\|_{\widetilde{L}^{1}_{T}(B^{\N}_{2,1})},\\
&\hspace{3,6cm}\lesssim \|\n q^{n-1}
\|_{\widetilde{L}^{\frac{4}{3}}_{T}(B^{\N+\frac{1}{2}}_{2,1})}\|\n q^{n-1}
\|_{\widetilde{L}^{4}_{T}(B^{\N-\frac{1}{2}}_{2,2})},\\
&\hspace{3,6cm}\lesssim \|q^{n-1}
\|_{\widetilde{L}^{\frac{4}{3}}_{T}(B^{\N+\frac{3}{2}}_{2,2})}\|q^{n-1}
\|_{\widetilde{L}^{4}_{T}(B^{\N+\frac{1}{2}}_{2,2})}.
\end{aligned}
$$
By following the same lines than the proof of theorem \ref{ftheo1}, we can easily conclude that $({\cal P}_{n})$ is verified. The only difficulty consists in treating the non linear term coming from the pressure. To deal with this term, we need to prove that $\rho^{n}$ belongs uniformly in $L^{\infty}$. From lemma \ref{petitlemme} and proposition \ref{flinear3}, we obtain that:
\begin{equation}
\begin{aligned}
&q^{n}_{L}\in\widetilde{L}^{\infty}(B^{\N}_{2,1})\cap\widetilde{L}^{1}(B^{\N+2}_{2,1}),\;{\rm div}u^{n}_{L}\in
\widetilde{L}^{\infty}(B^{\N-2}_{2,1})\cap\widetilde{L}^{1}(B^{\N}_{2,1})\\
&\hspace{6cm}\;\;\mbox{and}\;\;u^{n}_{L}\in
\widetilde{L}^{\infty}(B^{\N-1}_{2,1})\cap\widetilde{L}^{1}(B^{\N+1}_{2,1}).
\end{aligned}
\label{lineaireuL}
\end{equation}
It means in particular as $\bar{u}^{n}$ is in $\widetilde{L}^{1}(B^{\N+1}_{2,1})$,  and from (\ref{lineaireuL}) we conclude that ${\rm div}u^{n}$ is in $\widetilde{L}^{1}(B^{\N+1}_{2,1})$. It is then sufficient to prove that $q^{n}$ is in $L^{\infty}_{T_{n}}(L^{\infty})$. Indeed via the transport equation, we have:
$$\|q^{n}\|_{L^{\infty}}\leq \|q_{0}\|_{L^{\infty}}+\|{\rm div}u^{n}\|_{L^{1}(L^{\infty})}.$$
We can the deal with the pressure term as follows by using proposition \ref{composition}:
$$\|\frac{1}{\rho^{n}}\n P(\rho^{n})\|_{\widetilde{L}^{1}_{T}(B^{\N-1}_{2,1})}\leq C T\|q^{n}\|_{L_{T}^{\infty}(B^{\N}_{2,1})}.$$
We can now conclude easily  the proof of theorem \ref{ftheo2} by following the same lines as in theorem \ref{ftheo1}.\\
In the case of the existence of global strong solution with small initial data, we need to take in account the low frequencies. That is why we will include the pressure term in the linear part and we will use proposition \ref{1flinear3} to obtain estimate on $(q_{L},u_{L})$ where $(q_{L},u_{L})$ here stands for the  solution of:
\begin{equation}
\begin{cases}
\p_{t}q_{L}+{\rm div}u_{L}=0,\\
\p_{t}u_{L}-{\cal A}u_{L}+K\n q_{L}-\kappa\n(\D q_{L})=0,
\end{cases}
\label{lineaire}
\end{equation}
In particular it explains why we need of additional regularity on the initial density, i.e $q_{0}\in B^{\N-1}_{2,1}$. The rest of the proof for this case is similar to the case with large initial data.
  \hfill {$\Box$}
\section{Proof of theorem \ref{ftheo3}:}
We now want to prove theorem \ref{ftheo3} which is an ill-posedness theorem, in the sense that we have an explosion of the $L^{\infty}$ norm of the density in a arbitrary small time. To do this
we have just to choose a sequence of initial data $(q^{n}_{0},u_{0}^{n})_{n\in\mathbb{N}}$ such that:
\begin{itemize}
\item the sequel $(q^{n}_{0},\mathbb{P}u^{n}_{0},\mathbb{Q}u^{n}_{0})$ is uniformly bounded in
$(B^{\N-1}_{2,r}\cap B^{\N}_{2,r}\cap L^{\infty})\times B^{\N-1}_{2,2}\times  B^{\N-1}_{2,r}$ with $1\leq r\leq 2$,
\item we have also:
$$q^{n}_{0}\in B^{\N}_{2,1},\;\mathbb{P}u^{n}_{0}\in B^{\N-1}_{2,1}\;\;\mbox{and}\;\;\mathbb{Q}u^{n}_{0}\in B^{\N-1}_{2,1},$$
\item $(q^{n}_{0},\mathbb{P}u^{n}_{0},\mathbb{Q}u^{n}_{0})$ tends weakly to $(q_{0},\mathbb{P}u_{0},\mathbb{Q}u_{0})$ which is in $(B^{\N-1}_{2,r}\cap B^{\N}_{2,r}\cap L^{\infty})\times B^{\N-1}_{2,2}\times  B^{\N-1}_{2,r}$ ,
\end{itemize}
Furthermore we ask that there exists a sequel $(t_{n})_{n\in\mathbb{N}}$ such that $t_{n}\rightarrow_{n\rightarrow +\infty} 0$ and that the solution $(q^{n}_{L},u^{n}_{L})$ of system (\ref{lineaire}) with initial data $(q^{n}_{0},u_{0}^{n})$ verify:
$$\|{\rm div}u^{n}_{L}\|_{L^{1}_{t_{n}}(L^{\infty})}\rightarrow_{n\rightarrow +\infty} +\infty.$$
By theorem \ref{ftheo1}, we can show easily that it exists $T>0$ such that it exists strong solution $(q^{n},u^{n})$ verifying (\ref{impra}) with initial data $(q^{n}_{0},u_{0}^{n})$ on $(0,T)$. Furthermore by using the same arguments than theorem \ref{ftheo2}, we can show that:
$$q^{n}=q^{n}_{L}+\bar{q}^{n}\;\;\;\mbox{and}\;\;\;u^{n}=u^{n}_{L}+\bar{u}^{n},$$
such that $(\bar{q}^{n},\bar{u}^{n})$ is uniformly bounded in the following spaces:
\begin{equation}
\bar{q}^{n}\in\widetilde{L}_{T}^{\infty}(B^{\N}_{2,1})\cap\widetilde{L}_{T}^{1}(B^{\N+2}_{2,1}),\;\bar{u}^{n}\in\widetilde{L}^{\infty}(B^{\N-1}_{2,1})\cap\widetilde{L}^{1}(B^{\N+1}_{2,1}).
\label{explos}
\end{equation}
and such that by using proposition \ref{1flinear3} $(q_{L}^{n},u_{L}^{n})$ is uniformly bounded in the following spaces:
\begin{equation}
q_{L}^{n}\in\widetilde{L}^{\infty}(B^{\N}_{2,r})\cap\widetilde{L}^{1}(B^{\N+2}_{2,r}),u_{L}^{n}\in\widetilde{L}_{T}^{\infty}(B^{\N-1}_{2,r})\cap\widetilde{L}^{1}_{T}(B^{\N+1}_{2,2}).
\label{explos}
\end{equation}
We now want to consider the transport equation of system (\ref{NHV1}), we have then:
$$\p_{t}q_{n}=-{\rm div}u^{n}-u^{n}\cdot\n q^{n}.$$
We deduce that:
$$q_{n}(t,x)=-\int^{t}_{0}{\rm div}u^{n}(s,x)ds-\int^{t}_{0}u^{n}\cdot\n q^{n}(s,x)ds.$$
By paraproduct we can show that:
$$\|u^{n}\cdot\n q^{n}\|_{\widetilde{L}^{1}(B^{\N}_{2,1})}\leq \|q^{n}\|_{L^{\infty}_{T}(B^{\N}_{2,2})}\|u^{n}\|_{L^{1}_{T}(B^{\N+1}_{2,2})}.$$
We deduce that $u^{n}\cdot\n q^{n}$ is uniformly bounded in $\widetilde{L}^{1}(B^{\N}_{2,1})$ and then :
$$\big|\int^{t}_{0}u^{n}\cdot\n q^{n}(s,x)ds\big|\leq M.$$
As ${\rm div}u^{n}={\rm div}\bar{u}^{n}+{\rm div}u_{L}^{n}$ with ${\rm div}\bar{u}^{n}$ uniformly bounded in $ L^{1}_{T}(B^{\N}_{2,1})\h L^{1}_{T}(L^{\infty})$. Furthermore we recall that:
$$\|{\rm div}u^{n}_{L}\|_{L^{1}_{t_{n}}(L^{\infty})}\rightarrow_{n\rightarrow +\infty} +\infty.$$
It finally shows that:
$$\|q^{n}\|_{L^{\infty}_{t_{n}}(L^{\infty})}\rightarrow_{n\rightarrow +\infty} +\infty,$$
which is the desired result.
 \hfill {$\Box$}
\label{section5}
\section{Proof of theorem \ref{theo4}}
\label{section6}
For the proof of the existence of strong solution in finite time on $(0,T)$, it is a basic application of theorem 5 in \cite{fDD} except that we do not ask condition of smallness on the initial data. To do that we proceed as in \cite{aH2}. We just want to mention that we consider the momentum variable in order to choose some $p$ arbitrary big. Indeed it is crucial when we want to deal with the term ${\rm div}(u m)$ that the quantity $u m$ be inside of the divergence. It allows us to apply paraproduct laws.\\
We now want to prove the more interesting part of theorem \ref{theo4}, it means the blow-up criterion. The main idea is to obtain a gain of integrability on our solution $(\rho,v)$. By passing we can observe that $(\ln\rho, u)$ is also a solution of system (\ref{NHV1}) because we have a control on $\ln\rho$. After that we will use this gain of integrability on $v$ in order to obtain additional regularity on the density. It will be enough to conclude.\\
In the sequel we just will consider the case $N=3$, the case $N=2$ follows the same lines.
\subsubsection*{Gain of integrability on $v$} 
We multiply the the momentum equation of (\ref{3systeme2}) by $v|v|^{p-2}$ and integrate over $\R^{N}$, we obtain then:
\begin{equation}
\begin{aligned}
&\frac{1}{p}\int_{\R^{N}}\rho\p_{t}(|v|^{p})dx+\int_{\R^{N}}\rho u\cdot\n(\frac{|v|^{p}}{p})dx
+\int_{\R^{N}} \rho|v|^{p-2}|\n v|^{2}dx\\
&\hspace{1,5cm}+(p-2)\int_{\R^{N}} \rho\sum_{i,j,k}v_{j}v_{k}\p_{i}v_{j}\p_{i}v_{k}|v|^{p-4}dx+\int_{\R^{N}} |v|^{p-2}v\cdot\n\rho^{\gamma}dx=0.
\end{aligned}
\label{in1A}
\end{equation}
Next we observe that:
$$\sum_{i,j,k}v_{j}v_{k}\p_{i}v_{j}\p_{i}v_{k}=\sum_{i}(\sum_{j}v_{j}\p_{i}v_{j})^{2}=\sum_{i}\frac{1}{2}\p_{i}(|v|^{2}).$$
We get then as ${\rm div}(\rho u)=-\p_{t}\rho$ and by using (\ref{in1A}):
\begin{equation}
\begin{aligned}
&\frac{1}{p}\int_{\R^{N}}\p_{t}(\rho|v|^{p})dx+\int_{\R^{N}} \rho|v|^{p-2}|\n v|^{2}dx+(p-2)\int_{\R^{N}} \rho(\sum_{i}\p_{i}(|v|^{2})^{2}|v|^{p-4}dx\\
&\hspace{8cm}+\int_{\R^{N}} |v|^{p-2}v\cdot\n\rho^{\gamma}dx=0.
\end{aligned}
\label{in1A1}
\end{equation}
We have then by integrating over $(0,t)$ with $0<t<T$:
\begin{equation}
\begin{aligned}
&\frac{1}{p}\int_{\R^{N}}(\rho|v|^{p})(t,x)dx+\int^{t}_{0}\int_{\R^{N}} \rho|v|^{p-2}|\n v|^{2}(t,x)dtdx\\
&+(p-2)\int^{t}_{0} \int_{\R^{N}}\rho(\sum_{i}\p_{i}(|v|^{2})^{2}|v|^{p-4}(t,x)dtdx\leq \frac{1}{p}\int_{\R^{N}}(\rho_{0}|v_{0}|^{p})(x)dx\\
&\hspace{7cm}+|\int^{t}_{0} \int_{\R^{N}}|v|^{p-2}v\cdot\n\rho^{\gamma}(t,x)dtdx|.
\end{aligned}
\label{in1A1}
\end{equation}
By integration by part we have:
$$\int^{t}_{0} \int_{\R^{N}}|v|^{p-2}v\cdot\n\rho^{\gamma}(t,x)dtdx=-\int^{t}_{0} \int_{\R^{N}}{\rm div}(|v|^{p-2}v)a\rho(t,x)dtdx.$$
We have then:
$${\rm div}(|v|^{p-2}v)=|v|^{p-2}{\rm div}(v)+(p-2)|v|^{p-4}v\cdot(v\cdot\n v).$$
In the sequel we will note $p=N+\e$ with $\e>0$ and we will use the fact that as $(\frac{1}{\rho}-1)$ is in $L^{\infty}_{T}(\dot{B}^{0}_{N+\e,1}(\R^{N}))$, it is also in $L^{\infty}_{T}(L^{N+\e}(\R^{N})$. We can prove that $(\frac{1}{\sqrt{\rho}}-1)$ is in $L^{\infty}(L^{2p})$ because we have:
$$(\frac{1}{\sqrt{\rho}}-1)^{2p}=\big((\frac{1}{\rho}-1)+2(1-\frac{1}{\sqrt{\rho}})\big)^{p},$$
and:
$$
\begin{aligned}
\|\frac{1}{\sqrt{\rho}}-1\|_{L^{2p}}&\leq \|(\frac{1}{\rho}-1)+2(1-\frac{1}{\sqrt{\rho}})\|_{L^{p}}^{\frac{1}{2}},\\
&\leq C\big(1+\|(\frac{1}{\rho}-1)\|^{\frac{1}{2}}_{L^{p}}+\e\|(1-\frac{1}{\sqrt{\rho}})\|_{L^{2p}}+\|\frac{1}{\sqrt{\rho}}-1\|^{\frac{1}{2}}_{L^{1}}\big).
\end{aligned}
$$
We now recall that by (\ref{in1A1}), we control $\sqrt{\rho}\n(|u|^{\frac{p}{2}})$ in $L^{2}_{T}(L^{2})$ and as $(\frac{1}{\sqrt{\rho}}-1)$ belongs to $L^{\infty}(L^{2p})$ we can show that $\n(|u|^{\frac{p}{2}})$ is in $L^{2}_{t}(L^{q})$ with $\frac{1}{q}=\frac{1}{2}+\frac{1}{2p}$. \\
By Sobolev embedding we show easily that:
$$\||v|^{\frac{p}{2}-1}\|_{L^{\frac{2p}{p-2}}(L^{\frac{6p}{p-2}})}\leq C\|\n(|v|^{\frac{p}{2}})\|_{L^{2}(L^{2})}^{\frac{p-2}{p}},$$
and similarly that:
$$\||v|^{\frac{p}{2}-1}\|_{L^{\frac{2p}{p-2}}(L^{\frac{q_{1}p}{p-2}})}\leq C\|\n(|v|^{\frac{p}{2}})\|_{L^{2}(L^{q})}^{\frac{p-2}{p}}.$$
with $\frac{1}{q_{1}}=\frac{1}{2}+\frac{1}{2p}-\frac{1}{N}.$
By energy inequality (\ref{3inegaliteenergie1})  $\n\sqrt{\rho}$ belongs to $L^{\infty}_{T}(L^{2})$. By Sobolev embedding, we deduce as $(\sqrt{\rho}-1)$ is in $L^{\infty}(L^{1})$ that
$(\sqrt{\rho}-1)$ is in $L^{\infty}(L^{6})$.\\
\\
Finally, we have proves that $\sqrt{\rho}(\sqrt{\rho}|v|^{\frac{p}{2}-1}\n v)|v|^{\frac{p}{2}-1}$ is in $L^{1}(L^{\alpha}(\R^{N})$ with $\alpha>1$. We recall that:
$$\sqrt{\rho}(\sqrt{\rho}|v|^{\frac{p}{2}-1}\n v)|v|^{\frac{p}{2}-1}=(\sqrt{\rho}-1)(\sqrt{\rho}|v|^{\frac{p}{2}-1}\n v)|v|^{\frac{p}{2}-1}+(\sqrt{\rho}|v|^{\frac{p}{2}-1}\n v)|v|^{\frac{p}{2}-1}.$$
To do this we just apply H\"older's inequalities ,as by  (\ref{in1A1}) we have a control on $\sqrt{\rho}|v|^{\frac{p}{2}-1}\n v$ in $L^{2}_{T}(L^{2}(\R^{N}))$, on $(\sqrt{\rho}-1)$ in $L^{\infty}_{T}(L^{6}(\R^{N})$ and $|v|^{\frac{p}{2}-1}$ in $L^{\frac{2p}{p-2}}(L^{\frac{q_{1}p}{p-2}})$. We have then the following condition:
$$\frac{1}{\alpha}=\frac{p-2}{pq_{1}}+\frac{1}{2}+\frac{1}{6}=\frac{p-2}{p}(\frac{1}{2}+\frac{1}{2p}-\frac{1}{3})+\frac{1}{2}+\frac{1}{6}\leq 1,$$
because $p=3+\e>3$. By Sobolev embedding, we can exactly prove that $\sqrt{\rho}(\sqrt{\rho}|v|^{\frac{p}{2}-1}\n v)|v|^{\frac{p}{2}-1}$ is in fact in $L^{1}_{T}(L^{1})$ with:
$$
\begin{aligned}
&\big|\int^{t}_{0}\int_{\R^{N}}|v|^{p-2}v\cdot\n P(\rho)dxdt\big|\leq \|\n(|v|^{\frac{p}{2}})\|^{2-\e}_{L^{2}(L^{2})}(1+\|\frac{1}{\rho}-1\|_{L^{\infty}(L^{p})}\\
&\hspace{6cm}+\|\frac{1}{\sqrt{\rho}}-1\|_{L^{\infty}(L^{1})}+\|\sqrt{\rho}-1\|_{L^{\infty}(L^{1})})^{\beta},
\end{aligned}
$$
with $\beta$ big enough. By bootstrap in (\ref{in1A1}) we obtain finally that $\rho^{\frac{1}{p}}v$ is in $L^{\infty}(L^{p})$ for any $1\leq p<+\infty$.
\subsubsection*{Blow-up criterion}
We have then obtain that than for any $1\leq p<+\infty$, $\rho^{\frac{1}{p}}v$ belongs to $L^{\infty}(L^{p})$. We now want to consider the first equation of (\ref{3systeme2}):
\begin{equation}
\p_{t}q^{'}-\frac{\kappa}{\mu}\D q^{'}=-{\rm div}(\rho v),
\label{charegul}
\end{equation}
with $q^{'}=\rho-1$. We would like to take advantage of the gain of integrability that we have obtained on the velocity $v$ in order to show regularizing effect on the density. The goal is to transfer the information on the integrability of $v$ (which is a \textit{subscaling} estimate) on the density $\rho$. More precisely we have by proposition \ref{5chaleur} for any $1\leq p<+\infty$:
\begin{equation}
\|q^{'}\|_{\widetilde{L}_{T}^{\infty}(\dot{B}^{1}_{p,\infty})}\leq C(\|q^{'}_{0}\|_{\widetilde{L}^{\infty}(\dot{B}^{1}_{p,\infty})}+\|\rho v\|_{\widetilde{L}_{T}^{\infty}(\dot{B}^{0}_{p,\infty})}).
\label{incha}
\end{equation}
We know by energy inequalities (\ref{3inegaliteenergie1}) that:
\begin{equation}
(\rho-1)\in L^{\infty}(L^{\gamma}_{2})\;\;\;\mbox{and}\;\;\;\n\sqrt{\rho}\in L^{\infty}(L^{2}).
\end{equation}
We then split the product $\rho v$ as follows:
$$\rho v=(\rho^{1-\frac{1}{p}}-1)\rho^{\frac{1}{p}}v+\rho^{\frac{1}{p}}v=f(q^{'})\rho^{\frac{1}{p}}v+\rho^{\frac{1}{p}}v,$$
with $f(x)=(1+x)^{1-\frac{1}{p}}-1$. We now need to prove that $\rho v$ is in $\widetilde{L}_{T}^{\infty}(B^{0}_{p,\infty})$. As $\rho^{\frac{1}{p}}v$ is in $L^{\infty}_{T}(L^{p})$, we have then:
\begin{equation}
\begin{aligned}
\|f(q^{'})\rho^{\frac{1}{p}}v\|_{L_{T}^{\infty}(L^{p})}&\leq\|\rho^{\frac{1}{p}}v\|_{L_{T}^{\infty}(L^{p})}\|f(q^{'})\|_{L^{\infty}_{T}(L^{\infty})},\\
&\leq  C\|\rho^{\frac{1}{p}}v\|_{L_{T}^{\infty}(L^{p})}(1+\|q^{'}\|^{1-\frac{1}{p}}_{L^{\infty}_{T}(L^{\infty})}),\\
\end{aligned}
\label{prod1}
\end{equation}
But by interpolation as $p$ is big enough, we have for :
\begin{equation}
\|q^{'}\|_{L^{\infty}_{T}(L^{\infty})}\leq \|q^{'}\|^{1-\alpha}_{\widetilde{L}_{T}^{\infty}(\dot{B}^{1}_{p,\infty})}\|q^{'}\|^{\alpha}_{L^{\infty}_{T}(L^{2})},
\label{prod2}
\end{equation}
with $0<\alpha<1$. By injecting (\ref{prod1}) and (\ref{prod2}) in (\ref{incha}), we obtain:
\begin{equation}
\|q^{'}\|_{\widetilde{L}_{T}^{\infty}(\dot{B}^{1}_{p,\infty})}\leq C(\|q^{'}_{0}\|_{\widetilde{L}^{\infty}(\dot{B}^{1}_{p,\infty})}+\|\rho^{\frac{1}{p}} v\|_{L_{T}^{\infty}(L^{p})}\big(1+\|q^{'}\|^{1-\alpha}_{\widetilde{L}_{T}^{\infty}(\dot{B}^{1}_{p,\infty})}\|q^{'}\|^{\alpha}_{L^{\infty}_{T}(L^{2})}\big).
\label{incha1}
\end{equation}
As by (\ref{3inegaliteenergie1}) and Sobolev embedding, $\|q^{'}\|^{\alpha}_{L^{\infty}_{T}(L^{2})}$ is finite, we deduce from (\ref{incha1}) and from Young's inequalities that $q^{'}$ is in $\widetilde{L}_{T}^{\infty}(\dot{B}^{1}_{p,\infty})$ for any $N<p<+\infty$. It means that we have obtain that $\n\rho\in L^{\infty}_{T}(\dot{B}^{0}_{p,\infty})$. As in theorem \ref{theo4}, we assume that $\rho\in\widetilde{L}^{\infty}_{T}(\dot{B}^{0}_{N+\e,1})$ with $\e>0$, we have by paraproduct law (see proposition \ref{produit1}) when we are the critical case where the sum of the regularity index in nul:
$$
\|\n\ln\rho\|_{\widetilde{L}^{\infty}(\dot{B}^{0}_{p_{2},\infty})}\leq C(1+ \|\frac{1}{\rho}-1\|_{\widetilde{L}^{\infty}(\dot{B}^{0}_{N+\e,1})})\|\n\rho\|_{\widetilde{L}^{\infty}(B^{0}_{p,\infty})},$$
with $\frac{1}{p_{2}}=\frac{1}{p}+\frac{1}{N+\e}=\frac{1}{N+\frac{\e}{2}}$ (when we work with $p$ large enough).
We obtain then that $\n\ln\rho\in L^{\infty}(\dot{B}^{0}_{N+\frac{\e}{2},\infty})$ and that $\ln\rho\in L^{\infty}(\dot{B}^{1}_{N+\frac{\e}{2},\infty})$.\\
\\
Furthermore we have obtained that $\rho^{\frac{1}{p}}v$ belongs to $L^{\infty}(L^{p})$ for any $1\leq p<+\infty$ and that $\n\ln\rho\in L^{\infty}_{T}(\dot{B}^{0}_{N+\frac{\e}{2},\infty})$. As $u=v-\frac{\kappa}{\mu}\n\ln\rho$ and the fact that $(\frac{1}{\rho}-1)\in L^{\infty}_{T}(B^{0}_{N+\e})$, we obtain easily by following the same lines that for $p$ large enough  $u$ is in $L^{\infty}_{T}(\dot{B}^{0}_{N+\frac{\e}{2},\infty}$. To summarize what we have obtained, we have:
\begin{equation}
u\in L^{\infty}_{T}(\dot{B}^{0}_{N+\frac{\e}{2},\infty})\;\;\;\mbox{and}\;\;\;\ln\rho\in L^{\infty}_{T}(\dot{B}^{1}_{N+\frac{\e}{2},\infty}).
\label{hyperimp}
\end{equation}
We recall that the solution $(\rho, u)$ of system (\ref{3systeme2}) on $(0,T)$ are also solutions of system (\ref{NHV}) in the sense that $(q,u)$ is solution on $(0,T)$ of system (\ref{NHV}).\\
Easily we can prove that if $u_{0}\in \dot{B}^{0}_{N+\frac{\e}{2},\infty})$ and $\ln\rho_{0}=q_{0}\in \dot{B}^{1}_{N+\frac{\e}{2},\infty}$ then the system (\ref{NHV}) has a strong solution $(q,u)$ on $(0,T^{'})$ with:
$$T^{'}\geq \frac{C}{(1+\|u_{0}\|_{\dot{B}^{0}_{N+\frac{\e}{2},\infty}}+\|\ln\rho_{0}\|_{B^{1}_{N+\frac{\e}{2},\infty}})^{\beta}}.$$
This is an easy consequence of the fact that the initial data are choose subcritical, indeed in this case we have just to use proposition \ref{flinear3}, to prove that $(q_{L},u_{L})$ are small in function of the time $T^{«}$ in the critical norm.\\
It means that there exists a time $T^{'}\geq c>0$, where $c$ depends only on the physical coefficients and of subcritical initial data. We can construct by theorem \ref{ftheo1} a solution $(q_{1},u_{1})$ on $(T-\alpha,T-\alpha+T^{'})$ with initial data $(q(T-\alpha), u(T-\alpha))$ (here $\alpha<T^{'}$) which verifies (\ref{impra}) (in fact with a bit more regularity as we are subcritical) with Lebesgue index $p=N+\frac{\e}{2}$. The only difficulty is to prove that on $(T-\alpha,T)$ we have:
$$(q_{1},u_{1})=(q,u).$$
To do this, it suffices only to use the uniqueness part of theorem \ref{ftheo1}. It concludes the proof of theorem \ref{ftheo3}.
\hfill {$\Box$}
\section{Appendix}
\label{appendix}
In this appendix, we just give a technical lemma on the computation of the capillarity tensor.
\begin{lem}
When $\kappa(\rho)=\frac{\kappa}{\rho}$ with $\kappa>0$ then:
$${\rm div}K=\kappa\rho(\n\D\ln\rho+\frac{1}{2}\n(|\n\ln\rho|^{2})).$$
and:
$$
{\rm div}K=\kappa{\rm div}(\rho\n\n\ln\rho).
$$
\end{lem}
{\bf Proof:} We recall that:
$${\rm div}K
=\n\big(\rho\kappa(\rho)\D\rho+\frac{1}{2}(\kappa(\rho)+\rho\kappa^{'}(\rho))|\n\rho|^{2}\big)
-{\rm div}\big(\kappa(\rho)\n\rho\otimes\n\rho\big).$$
When $\kappa(\rho)=\frac{\kappa}{\rho}$, we have:
\begin{equation}
{\rm div}K=\kappa\n\D\rho-\kappa{\rm div}(\frac{1}{\rho}\n\rho\otimes\n\rho).
\label{1a1}
\end{equation}
But as:
$$\D\rho=\rho\D\ln\rho+\frac{1}{\rho}|\n\rho|^{2},$$
we have by injecting this expression in (\ref{1a1}):
\begin{equation}
{\rm div}K=\kappa\rho\n\D\ln\rho+\kappa\n\rho\D\ln\rho+\kappa\n(\frac{1}{\rho}|\n\rho|^{2})-\kappa{\rm div}(\frac{1}{\rho}\n\rho\otimes\n\rho).
\label{a2}
\end{equation}
As we have:
$$\kappa{\rm div}(\frac{1}{\rho}\n\rho\otimes\n\rho)=\kappa\D\ln\rho\n\rho+\n(\frac{1}{\rho}|\n\ln\rho|^{2})-\frac{\kappa}{2}\rho\n(|\n\ln\rho|^{2}).$$
It concludes the first part of the lemma.\\
\\
We now want to prove that we can rewrite (\ref{cap}) under the form of a viscosity tensor. To see this, we have:
$$
\begin{aligned}
{\rm div}(\rho\n(\n \ln\rho))_{j}&=\sum_{i}\p_{i}(\rho\p_{ij}\ln\rho),\\
&=\sum_{i}[\p_{i}\rho\p_{ij}\ln\rho+\rho\p_{iij}\ln\rho],\\
&=\rho(\D\n \ln\rho)_{j}+\sum_{i}\rho\p_{i}\ln\rho\p_{j}\p_{i}\ln\rho),\\
&=\rho(\D\n \ln\rho)_{j}+\frac{\rho}{2}(\n(|\n\ln\rho|^{2}))_{j},\\
&={\rm div}K.
\end{aligned}
$$
We have then:
$${\rm div}K=\kappa{\rm div}(\rho\n\n\ln\rho)=\kappa{\rm div}(\rho D(\n\ln\rho)).$$
\hfill {$\Box$}

\end{document}